\documentclass[aos,MSNbibl,seceqn,dvips]{arximspdf}
\usepackage{mathbh}

\doi{10.1214/12-AOS1055} 
\volume{40}
\issue{6}
\pubyear{2012}
\firstpage{2910}
\lastpage{2942}

\makeatletter

\newcommand{\rright}{\right}
\newcommand{\lleft}{\left}
\newcommand{\rrvert}{\vert}
\newcommand{\llvert}{\vert}

\newtheorem{theorem}{Theorem}[section]
\newtheorem{corollary}[theorem]{Corollary}

\def\definedby{\stackrel{\Delta}{=}}

\def\indic{{\mathbh{1}}}
\def\a{\alpha}
\def\cal{\mathcal}

\def\P{{\mathbb P}}
\def\E{{\mathbb E}}
\def\Ee{{\mathbb E}}
\def\real{{\mathbb{R}}}
\def\definedas{\stackrel{\Delta}{=}}
\def\Min{{\cal M}}
\def\implies{\Rightarrow}

\newcommand{\vari}{\operatorname{Var}}
\newcommand{\Tr}{\operatorname{Tr}}
\newcommand{\mink}{{\cal M}}
\newcommand{\hauss}{{\cal H}}
\newcommand{\lips}{{\cal L}}
\newcommand{\detr}{\operatorname{detr}}

\makeatother

\begin{document}
\begin{frontmatter}

\title{Rotation and scale space random fields and the
Gaussian kinematic formula}
\runtitle{Rotation and scale space random fields}

\begin{aug}
\author[A]{\fnms{Robert J.} \snm{Adler}\thanksref{t1,t5,t4}\ead[label=e1]{robert@ieadler.technion.ac.il}\ead[label=u1,url]{webee.technion.ac.il/people/adler}},
\author[A]{\fnms{Eliran} \snm{Subag}\corref{}\thanksref{t5}\ead[label=e2]{elirans@techunix.technion.ac.il}}
\and\\
\author[B]{\fnms{Jonathan E.} \snm{Taylor}\thanksref{t1,t3,t4}\ead[label=e3]{jonathan.taylor@stanford.edu}\ead[label=u3,url]{www-stat.stanford.edu/\textasciitilde jtaylo/}}
\runauthor{R. J. Adler, E. Subag and J. E. Taylor}
\affiliation{Technion, Technion and Stanford University}
\address[A]{R. J. Adler \\
E. Subag\\
Faculty of Electrical Engineering\\
Technion, Haifa 32000\\
Israel \\
\printead{e1}\\
\hphantom{E-mail: }\printead*{e2}\\
\printead{u1}}
\address[B]{J. E. Taylor\\
Department of Statistics\\
Stanford University\\
Stanford, California 94305-4065\\
USA \\
\printead{e3}\\
\printead{u3}} 
\end{aug}

\thankstext{t1}{Supported in part by US--Israel Binational Science Foundation, 2008262.}

\thankstext{t5}{Supported in part by Israel Science Foundation, 853/10, and AFOSR FA8655-11-1-3039.}

\thankstext{t3}{Supported in part by NSF Grant DMS-09-06801.}

\thankstext{t4}{Supported in part by NSF-SGER Grant DMS-08-52227.}

\received{\smonth{12} \syear{2011}}
\revised{\smonth{6} \syear{2012}}

%
\begin{abstract}
We provide a new approach, along with extensions, to results in
two important papers of Worsley, Siegmund and coworkers
closely tied to the statistical analysis of fMRI (functional magnetic
resonance imaging) brain data. These papers studied approximations for
the exceedence probabilities of scale and rotation space random fields,
the latter playing an important role in the statistical analysis of
fMRI data.
The techniques used there came either from the Euler characteristic
heuristic or via tube formulae, and to a large extent were carefully
attuned to the specific examples of the paper.

This paper treats the same problem, but via calculations based on the
so-called Gaussian kinematic formula. This allows for extensions of the
Worsley--Siegmund results to a wide class of non-Gaussian cases. In addition,
it allows one to obtain results for rotation space random fields in any
dimension via reasonably straightforward
Riemannian geometric calculations. Previously only
the two-dimensional case could be covered, and then only via computer
algebra.

By adopting this more structured approach to this particular problem,
a~solution path for other, related problems becomes clearer.
\end{abstract}

%

\end{frontmatter}

\section{Introduction}\label{secIntroduction}

In this paper we provide a new approach, along with extensions, to
results in
two important papers closely tied to the statistical analysis of fMRI
(functional magnetic
resonance imaging) brain data. The first, by David Siegmund and the late
Keith Worsley~\cite{S-WOR95}, treated the problem of testing for a
spatial signal\vspace*{1pt}
with unknown location and scale perturbed by
a Gaussian random field on $\mathbb{R}^N$.
While their motivation came from fMRI, their results are also relevant to
areas such as astronomy and genetics, or indeed to many applications in which
one searches a large space of noisy data
for a relatively small number
of signals which manifest as ``bumps'' in the random noise.

In a second paper, together with Shafie and Sigal
\cite{Shafie03rotationspace}, they treated an extended version of the same
basic problem, allowing the asymmetric bumps to also have a direction
associated to them, in addition to scale changes that were
direction dependent. This added an additional ``rotation parameter'' to
the model and
complicated it to the extent that they were unable to treat anything
beyond the
two-dimensional ($N=2$) case. Thus, for example, while their results
could be used
to examine cortical surface images, they could not be applied to
three-dimensional
fMRI analysis.

For etymological reasons,
the
first of these cases was called the \textit{scale space} case and the second
the \textit{rotation space} case. In both cases, the mathematics of these papers
concentrated on approximating the exceedence probabilities
%
\begin{equation}\label{exceedence1equn} 
\P \Bigl\{\sup_{t\in M} f(t) \geq u \Bigr\}
\end{equation}
for large $u$, a smooth, Gaussian, real-valued, random field $f$, and some
parameter space $M$. For reasons that will become clearer
soon, in the $N$-dimensional scale space case, the effective dimension
of $M$ was $N+1$, while in the rotation space case it jumped to
$N+N(N+1)/2$ (or 5, for the
only case they could handle, when $N=2$).

The aims of the current paper are twofold. First, we shall generalize
the results of the
rotation space to any dimension, but by a very different method.
Rather than taking a ``from first principles''
approach, as was the case in~\cite{Shafie03rotationspace}, we shall
compute everything
from a basic result known as the \textit{Gaussian kinematic formula}
[GKF, cf. (\ref{pmainequation}) below],
using techniques of Riemannian geometry. In fact, we shall derive the rotation
space results from a more general setting, which is of interest in itself.
This is quite different from the analysis in
\cite{Shafie03rotationspace}, in which calculations were so involved
that they could only be handled
by symbolic calculus using computers and, then, as already mentioned,
only in two dimensions. According to the ``no free lunch'' principle,
there is a price to pay
here, in that the specific
$N=2$ approximations of~\cite{Shafie03rotationspace} have a
slightly higher level of precision.

Second, we shall recover the results of the scale space case, again
using the GKF.
We shall do this in two ways. First, we shall obtain it as an immediate
corollary to the
``more general setting'' mentioned above. This is quick, but once again
does not give
the best results. However, we shall also derive the ``best'' results from
a more detailed calculation, something which seems too difficult to do for
rotation space fields.

At first this might seem as if we are attempting to kill the proverbial
ant with a
sledge hammer, but it turns out that there are at least three
advantages to this approach. First of all, the ant in question is a
rather large
one, and any way of dealing with it requires considerable effort. The
basic calculations in both
\cite{S-WOR95} and here take about 10 pages, and the
essential difference
is whether one prefers to work with Riemannian topology or integral geometry,
hard calculus and occasionally symbolic (machine) calculus. This may be
more a question of mathematical taste than anything else. The second
advantage of using the
GKF as the basic starting point, which is significant and no longer
merely an issue of taste, is that while
the results of~\cite{Shafie03rotationspace} and
\cite{S-WOR95} apply,
as stated, only to Gaussian fields, the GKF allows immediate extension
to a much
wider class of examples. Section
\ref{gkfsubsec} explains this.

Finally, both of these cases are instructive in terms of how to use the
GKF in
more general settings,
and our hope is that presenting them in detail will provide examples
which will make
the GKF easier to use in the future for other, perhaps quite different,
examples.

We shall say no more in the way of motivation for this paper.
Discussions and
examples can be found in both~\cite{Shafie03rotationspace} and
\cite{S-WOR95} as well as the papers cited there, and the many
more recent papers
citing~\cite{Shafie03rotationspace} and~\cite{S-WOR95}.
In particular, it is explained there why (\ref{exceedence1equn}) is
an important
quantity in statistical hypothesis testing.
Rather, we shall, in the following section, immediately turn to formal
definitions
of scale and rotation random fields, a description of the GKF and its relation
to~\cite{Shafie03rotationspace} and~\cite{S-WOR95}, and a
statement of the main results of the paper in Section~\ref{lkcrotatesubsec}.
Fuller results for scale space fields, which go considerably
beyond those of Section~\ref{lkcrotatesubsec}, are given at the end
of the paper,
in Section~\ref{scalesec}.

Understanding or using the results of the paper does
not require a deeper knowledge of differential
geometry beyond that which
we describe along the way.
This is not true of the proofs, and in Section~\ref{beforesection}
we point to what else is needed. The remainder of the paper is
taken up with proofs.

\section{Random fields, the GKF and our main results}
\label{mainsec}

The common element behind all the random fields in this paper is
Gaussian white noise.

\subsection{Gaussian white noise}
\label{noisesubsec}
The simplest way to define Gaussian white noise $W$ on $\mathbb{R}^N$
is as a random
field indexed by the Borel sets $\cal B^N$, finitely additive, taking
independent values
on disjoint sets, and such that, for all $A\in\cal B^N$ with finite
Lebesgue measure
$\lambda_N(A)\equiv|A|$,
$ W(A) \sim N(0, |A|)$. %
It is standard fare to define integrals
$\int_{\mathbb{R}^N} \varphi(t) W(dt)$,
for \textit{deterministic} $\varphi\in L_2(\mathbb{R}^N)$,
which are mean zero,
Gaussian random variable with covariances
%
\begin{equation}
\label{introinnerprodallequation} \E \biggl\{\int_{\mathbb{R}^N}
\varphi(t) W(dt) \int_{\mathbb
{R}^N} \psi(t) W(dt) \biggr\} =  \int
_{\mathbb{R}^N} \varphi(t)\psi(t) \lambda_N(dt).
\end{equation}

\subsection{Scale space fields}
\label{scalespacesubsec}

The Gaussian scale space random field is obtained by smoothing white
noise with
a spatial filter over a range of filter widths or scales. Formally, for
a ``filter'' $h\in L_2(\mathbb{R}^N)$ it is defined as
%
\begin{equation}
\label{gsst} f(\sigma,t)=\sigma^{-N/2} \int_{\real^N} h
\biggl(\frac{t-u}{\sigma
} \biggr) \,dW(u).
\end{equation}
Typically $t$ ranges over a nice subset $T$ of $\mathbb{R}^N$, and we take
$\sigma\in[\sigma_l,\sigma_u]$, for finite $0<\sigma_l<\sigma_u$.

The parameter
space therefore includes both scale and location parameters \mbox{\cite
{S-WOR95,WOR98}}.
From a statistical point of view,
the scale space field is a continuous wavelet transform of white noise
that is designed to be powerful at detecting a localized signal of unknown
spatial scale and location.

This concludes the definition, and we now turn to some additional
assumptions, as
well as setting up some notation and further definitions that we shall
require to
state the main results and to carry out the proofs.

For convenience, and without loss of generality,
we shall assume that the filter $h$ is normalized so that
$\int_{\mathbb{R}^N} h^2(u) \,du= 1$, which will ensure that the
variance of
$f$ is also one.
We shall also assume that $h$ is twice continuously differentiable on
$\mathbb{R}^N$
and that all partial derivatives are in $L_2(\mathbb{R}^N)$.

With some loss of generality, that, in practice, is not generally too
restrictive,
we shall assume that there exists a $\gamma>0$ such that
%
\begin{equation}
\label{ECscaledequn} \int_{\mathbb{R}^N} \nabla h(u) \bigl(\nabla h
(u)\bigr)' \,du= \gamma I_{N\times N},
\end{equation}
where we write $I$ or $I_{N\times N}$ to denote the $N\times N$
identity matrix.
(Note: We take all our vectors to be column vectors.)
Equation (\ref{ECscaledequn}) will follow, for example, if the
filter is spherically symmetric.
Two common examples are given by the standard Gaussian kernel
$
h(t) = {\pi^{-N/4}}e^{-|t|^2/2}$
and the Marr or ``Mexican hat'' wavelet
\[
h(t)= \biggl[\frac{4N}{(N+2)\pi^{N/2}} \biggr]^{1/2} \biggl(1-\frac
{|t|^2}{N}
\biggr) e^{-|t|^2/2}.
\]
It is easy to check that in the Gaussian case (\ref{ECscaledequn})
holds with
$\gamma=1/2$ and in the Marr wavelet
case with $\gamma=(N+4)/(2N)$. Furthermore,
the covariance function of the scale space field
is easily calculated, via (\ref{introinnerprodallequation}),
to be
\[
C\bigl((\sigma_1,t_1),(\sigma_2,t_2)
\bigr)= 
\frac{1}{(\sigma_1\sigma_2)^{N/2}}\int
_{\mathbb{R}^N} h \biggl(\frac
{t_1-t_2-v}{\sigma_1} \biggr) h \biggl(
\frac{-v}{\sigma_2} \biggr) \,dv.
\]
Note that, for fixed $\sigma$, it follows from the previous line that
$f$ is stationary in $t$. However, $f$ is definitely not stationary as
a process
in the pair $(\sigma,t)$.

To apply the GKF and to state our results, we shall
need information about the derivatives
of $f$, which exist due to our assumptions on $h$. Either from first principles,
by differentiating the covariance function
or by using (5.5.5) of~\cite{TheBook}, we find
that all first-order partial derivatives of $f$ with respect to
the space variables are uncorrelated with the first-order derivative in
the scale
variable. This lack of correlation, which of course is equivalent to
independence since
all fields are Gaussian, is \textit{crucial} to the proofs, and without
it, it is unlikely that
we could complete the detailed calculations that we shall carry out
later. The same techniques also give
%
\begin{eqnarray}
\label{ECscale-variancesequn} \kappa&\definedby& {\sigma}^2 \cdot
\vari \biggl(\frac{\partial
f(\sigma,t)}{\partial\sigma} \biggr) = \int_{\mathbb{R}^N} \bigl[
\bigl\langle u,\nabla h (u)\bigr\rangle+ N h(u)/2 \bigr]^2 \,du,
\\
\label{ECscale-variances2equn}
\Lambda_{\sigma} &\definedby& \vari \bigl( \nabla_t f(
\sigma,t) \bigr) = \sigma^2 \int_{\mathbb{R}^N} \nabla h(u)
\bigl(\nabla h(u)\bigr)^\prime \,du,
\end{eqnarray}
where $\nabla_t f$ denotes the gradient of $f$ with respect to the
elements of $t$ only.
In view of (\ref{ECscaledequn}), we have $\Lambda_{\sigma}=\sigma^2 \gamma I_{N\times N}$. We should note that it is not obvious at this
point that
the quantity $\sigma^2 \vari (\partial f(\sigma,t)/\partial
\sigma )$ should be independent of $\sigma$. However, this is
in fact a consequence
of calculations in Section~\ref{secscale}. For the Gaussian and Marr
wavelet kernels the values of $\kappa$ are $1/2$ and $N/2$, respectively.

Excluding Section~\ref{scalesec}, from now on we shall assume that
$h$ is
spherically symmetric, so that we can define
$k\dvtx\real_+\to\real$ by
%
\begin{equation}
\label{eqkernel} h(x)  = k\bigl(\|x\|^2\bigr).
\end{equation}
In this case, $\gamma$, $\kappa$ and other constants which will
appear throughout can be expressed in terms of $k$. In particular,
retaining consistency
with
(\ref{ECscaledequn}) and (\ref{ECscale-variancesequn})
and introducing three new terms, we have
%
\begin{eqnarray}
\label{eqconstants}
\gamma &=& 4 \int_{\mathbb{R}^{N}}
\dot{k}^{2} \bigl( u'u \bigr) u_{1}^{2}
\,du,\qquad \kappa_{2,2}=\int_{\mathbb{R}^{N}}\dot{k}^{2}
\bigl( u'u \bigr) u_{1}^{2}u_{2}^{2}
\,du,
\nonumber\\[-8pt]\\[-8pt]
\kappa_{4}&=&\int_{\mathbb{R}^{N}}\dot{k}^{2}
\bigl( u'u \bigr) u_{1}^{4} \,du, \qquad\rho=\int
_{\mathbb{R}^{N}}k \bigl( u'u \bigr) \dot{k} \bigl(
u'u \bigr) u_{1}^{2} \,du
\nonumber
\end{eqnarray}
and also define
%
\begin{equation}
\label{eqdefCtildek} \widetilde{C}_{k,N}^2 \definedby
\frac{\kappa}{4} = N^2 \biggl(\frac{1}{16} +
\frac{\rho}{2} + \kappa_{2,2} \biggr) + 2N\kappa_{2,2}.
\end{equation}

Note for later use that
%
\begin{equation}
\label{eqkappaidentity} \kappa_{4}-3\kappa_{2,2}=\int
_{\mathbb{R}^{N}}\dot{k}^{2} \bigl( u'u \bigr)
\bigl(u_{1}^{4}-3u_{1}^{2}u_{2}^{2}
\bigr) \,du = 0.
\end{equation}
This can be easily proved for a general isotropic kernel [instead of
$\dot{k}^{2} ( u'u )$ above] by using
Fubini's theorem to restrict to the two-dimensional case and then
expressing the integral using spherical coordinates.

\subsection{Rotation space fields}
\label{rotationspacesubsec}
As already described, rotation space random fields are based on the
same underlying white noise
model as scale space fields, but also allow for rotations and scale
changes that are
direction dependent. More precisely, write
\[
\mathbb{S}=\mathbb{S}\bigl(\sigma_2^{-2},
\sigma_1^{-2}\bigr)\triangleq \bigl\{ S\in
\mathrm{Sym}_{N}^+\dvtx\sigma_{2}^{-2}I\leq S\leq
\sigma_1^{-2}I \bigr\},
\]
where $\mathrm{Sym}_N^+ \subset\mathrm{Sym}_{N}$ is the set of
$N\times
N$ positive
definite symmetric matrices, $0<\sigma_1<\sigma_2<\infty$,
$K\definedas N(N+1)/2$, $D\definedas N+K$ and for two matrices $A$ and
$B$ we write $A \leq B$ if, and only if, $B-A$ is nonnegative definite.
Then the rotation space random field, defined over the $D$-dimensional space
$M\definedas T\times\mathbb{S}$, is given by
%
\begin{equation}
\label{rsfdefn} f (t,S )=\llvert S\rrvert^{{1}/{4}}\int
_{\mathbb
{R}^{N}}h \bigl( S^{{1}/{2}} (t-u ) \bigr) \,dW (u ).
\end{equation}
Again, we shall assume that
$h(t)$ is spherically symmetric, twice continuously differentiable on
$\mathbb{R}^N$,
that all partial derivatives are in $L_2(\mathbb{R}^N$)
and normalized so that $\int_{\mathbb{R}^N} h^2(u) \,du= 1$. Note
that the
eigenvalues and eigenvectors of $S$ determine the scalings and
the corresponding ``principal directions'' for the
rotation space field.

The constant analogous to $\widetilde{C}_{k,N}$ in rotation
space is
%
\begin{equation}
\label{eqdefCk} %
C_{k,N} \definedby (2\pi
)^{{D}/{2}-N}\prod_{j=1}^{N}
\frac{\Gamma (1/2 )}{\Gamma (j/2 )}c_{k,N},
\end{equation}
where
%
\begin{equation}
\label{eqdefck} %
c_{k,N}^2
\definedby (2\kappa_{2,2} )^{K} \biggl(1+N (2
\kappa_{2,2} )^{-1} \biggl(\frac{1}{16}+
\frac{\rho
}{2}+\kappa_{2,2} \biggr) \biggr).
\end{equation}

Reassuringly, when $N=1$, in which case scale space is equivalent to rotation
space, it is easy to check that $C_k=\widetilde{C}_k$. Another
constant which will appear later is
%
\begin{equation}
\label{eqdefDk}
D_{k,N}^2
\definedby 
2\kappa_{2,2}
\frac{\widetilde{C}_{k,N}^2/N}{\widetilde{C}_{k,N-1}^2/(N-1)}.
\end{equation}

As for scale space fields, we shall once again need information about
variances and
covariances of the first-order derivatives of rotation space fields.
This time, however,
the formulae are somewhat more complicated, and, since one of the
parameters of
$f(t,S)$ is a matrix, we shall have to explain what we mean by
differentiation with
respect to this parameter. Details will be given below in Section \ref
{dependencesubsec}.
At this point we note only that, as before, derivatives with respect to
$t$ and $S$ are
independent, and, again, this is crucial to our ability to carry out detailed
computations. Furthermore, although once again $f$ is neither
stationary nor isotropic,
it does have zero mean and constant variance.

\subsection{Regularity conditions}
\label{regularitysubsec}

Before we state the GKF, we need to say something about the regularity
conditions
needed for it to hold. These are of three kinds, relating to the
smoothness of the
scale and rotation space fields, the structure of $T$ (and so of $M$)
and the smoothness of a
transformation $F$ that we shall meet only in the next subsection. Full
details are given in Chapters 11 and 12 of
\cite{TheBook}.

We start with $T$, which we take to be a compact domain in $\mathbb
{R}^N$ with
a $C^2$ boundary (although
\cite{TheBook} would allow us to assume less).
For the remainder of this paper we require no convexity conditions on
$T$, although
in order to apply the Euler characteristic heuristic of
Section~\ref{motivationsec}
it is necessary to assume that $T$ is locally convex in the
sense of Definition 8.2.1 of~\cite{TheBook}.

The second set of conditions is required to ensure that $f$ is, with
probability one, $C^2$,
as a function of each
of the space and scale or rotation variables, and also nondegenerate in
a certain
(Morse theoretic) sense. Sufficient conditions under which this holds
are given in Sections
11.3 and 12.1 of~\cite{TheBook}.
It is easy to check that these conditions are satisfied for the
Gaussian, Marr
and other smooth kernels.

\subsection{The Gaussian kinematic formula}
\label{gkfsubsec}
We shall state the GKF in more generality than needed for this paper,
since to do so requires
little more than an additional sentence or two.
Suppose $M$ is a $C^2$, Whitney
stratified manifold satisfying mild side conditions
and $D$ a similarly nice stratified submanifold of $\mathbb{R}^k$.
For the exact definitions of ``mild'' and ``nice'' see~\cite{TheBook},
but, for the purposes
of this paper, it suffices that $M$ is either a scale or rotation
space, with $T$ satisfying the
conditions of the previous subsection. Also, if $F\dvtx\mathbb
{R}^k\to\real$
is $C^2$, then it will suffice
that $D=F^{-1}([u,\infty))$ for some~$u$.

Let
$f=(f^1,\ldots,f^k)\dvtx M\to\real^k$ be
a vector valued random process, the components of which are independent,
identically distributed, real-valued, $C^2$, centered,
unit variance, Gaussian processes, satisfying the conditions of the
previous subsection.
The Gaussian kinematic formula (due
originally to Taylor in~\cite{Taylor-thesis,Taylor1} and extended in
\cite{TheBook,TAannals})
states that, for certain additive set functionals $\lips_0,\ldots,\lips_{\dim M}$,
%
\begin{equation}
\label{pmainequation}\qquad
\E \bigl\{\lips_i \bigl(M\cap
f^{-1}(D) \bigr) \bigr\}= \sum_{j=0}^{\dim M -i}
\lleft[\matrix{i+j
\cr
j} \rright] (2\pi )^{-j/2}\lips_{i+j}(M)
\Min_j^{k}(D).
\end{equation}
There is a lot to explain here. The ``flag''
coefficients
$ \bigl[{n \atop j} \bigr]
= {n\choose j} {\omega_n}/{\omega_{n-j} \omega_j}$,
where $\omega_n$ is the volume of the
unit ball in $\real^n$. The $\Min_j^{k}(D)$, known as the Gaussian
Minkowski functionals of $D$, are determined via the tube expansion
%
\begin{equation}
\label{tubesgaussexp1equn1} \P \bigl\{ \xi\in\bigl\{y\in\mathbb{R}^k
\dvtx\operatorname {dist}(y,D)\leq\rho\bigr\} \bigr\} = \sum
_{j=0}^{\infty} \frac{\rho^j}{j!} \Min_j^{k}(D),
\end{equation}
where $\xi\sim N(0,I_{k\times k})$ and dist is the usual Euclidean distance.
In the case $D=F^{-1}([u,\infty))$ it typically involves no\vadjust{\goodbreak} more than
multivariate
calculus to compute the $\Min_j^{k}(D)$, and many examples are given
in~\cite{TheBook}
and~\cite{ARF}. We shall give one example in a moment.

The most important elements of (\ref{pmainequation}) are the $\lips_j$, known as
Lipschitz--Killing curvatures (LKCs).
The first of these is rather simple, as it is always the Euler
characteristic $\chi$, that is,
$
\lips_0 (A) \equiv\chi(A)$
for all nice $A$. Thus, (\ref{pmainequation}), with $i=0$, provides
us with
an expression for
$\E\{\chi (M\cap f^{-1}(D) )\}$, which is what is needed
to approximate exceedance
probabilities, as explained in the following section.

For a simple but extremely important example, suppose
that $k=1$ and $F$ is the identity.
Then it is
easy to compute the $\mink_j^1([u,\infty))$ to see that, with $i=0$,
and $f\dvtx M\to\real$ satisfying the conditions above, (\ref{pmainequation})
becomes
%
\begin{equation}
\label{eecgaussianequn} E \bigl\{ \chi \bigl(A_{u}(f,M) \bigr) \bigr
\} =e^{-
{u^{2}}/{2}}\sum_{j=0}^{\dim M} (2\pi
)^{-
({j+1})/{2}}\mathcal{L}_{j} (M )H_{j-1} (x ),
\end{equation}
where
$A_u(f,M)= \{t\in M\dvtx f(t)\geq u\}$
is an excursion set,
the $H_{n}$, $n\geq0$, are Hermite polynomials,
and $H_{-1}(x) \definedas[\varphi(x)]^{-1}\int_x^\infty\varphi(u)
\,du$, where
$\varphi$ is the standard Gaussian density.
Since the structure of $M$ affects only the LKCs, and neither the Minkowski
functionals nor the flag coefficients, computing LKCs is an independent
problem. This is
what this paper is about, for scale and rotation space fields.

(A word to the purely applied statistician who needs only
a numerical, data driven method for estimating LKCs from data and
wishes to avoid
theory. Two different
such methods~\cite{Adler-Bartz-Kou,Taylor-Worsley-JASA} are currently
available.)

Throughout, we assume that $M$ is a regular (in the sense of \cite
{TheBook}) stratified manifold.
These are basically sets that can be partitioned into a
disjoint union of $C^2$ manifolds, so that we can write
%
\begin{equation}
\label{ppartitionequation} M = \bigsqcup_{j=0}^{\dim M}
\partial_j M,
\end{equation}
where each nonempty stratum, $\partial_jM$, $0\leq j\leq\dim(M)$, is a
$j$-dimensional manifold whose closure contains $\partial_iM$ for all $i<j$.
For a $N$-dimensional cube, $\partial_NM$ is its
interior, $\partial_{n-1}M$ is the union of its ($N-1$)-dimensional
(open) faces and so
on, down to $\partial_0M$, which is the collection of its vertices.

If such a space is locally convex, then a simple way to define a
Euclidean version $\cal L^E_j$ of
the LKCs is via a \textit{tube formula},
which states that, for small enough $\rho$, and $M \subset\real^N$,
%
\begin{equation}
\label{geometrysteinerformula} \cal H^E_N \Bigl( \Bigl
\{t\in\real^N\dvtx\inf_{s\in M} \|t-s\| \leq\rho \Bigr\} \Bigr)
= \sum_{j=0}^N \omega_{N-j}
\rho^{N-j} \lips_j^E(M),
\end{equation}
where
(for reasons that will become clearer later)
we write $\cal H^E_N$ for the Lebesgue (Hausdorff) measure on $\mathbb{R}^N$.
It follows that, for rectangles,\break
$\lips_j^E(\prod_{i=1}^{N} [0,T_i])
= \sum T_{i_1}\cdots T_{i_j},
$
where the sum is taken over all 
distinct choices of
subscripts.

Note that it will always be true that $\lips_N^E(M)$ is the volume of
$M$, while
$\lips_{N-1}^E(M)$ is half the surface measure of $M$. Furthermore,
$\lips^E_0(M)=\chi(M)$.

In what follows, however, we shall need to work in a non-Euclidean
setting, as the
Lipschitz--Killing curvatures in the GKF, and thus in (\ref
{eecgaussianequn}), are computed relative to a Riemannian
metric induced by a Gaussian field. This metric $g$ on $M$ is defined by
%
\begin{equation}
\label{metric} g_x(X_x,Y_x)
\definedas \E\bigl\{\bigl(X_xf^i_x\bigr)
\bigl( Y_xf^i_x\bigr)\bigr\}
\end{equation}
for any $i$ and for $X_x,Y_x\in T_xM$, the tangent space to $M$ at
$x\in M$.
This can also be written in terms of the covariance function as
%
\begin{equation}\label{eqRiemannianmetric}
g_{x} (X_{x},Y_{x} )=X_{x}Y_{y}C
(x,y )|_{x=y}.
\end{equation}

Both of the Gaussian fields of interest to us---viz. the scale and rotation
space fields---which we defined on compact sets, can easily be
extended, to
$\mathbb{R}^N\times\real^+$ and $\real^N \times\mathrm{Sym}_{N}^+$,
respectively, and so
(\ref{metric}) defines metrics on these spaces.
In all that follows, unless stated otherwise,
these metrics, or their restrictions, are what we shall work with, and
it will
be in relation to them that Hausdorff measures
$\mathcal{H}_j$ and Lipschitz--Killing curvatures $\lips_j$ are
defined. When we
consider a submanifold $L$ of either of these spaces the notation
$\mathcal{H}_j^L$, $\lips_j^L$ will be used. When the standard
Euclidean metric is
used these will be replaced by $\mathcal{H}_j^E$ and $\lips_j^E$.

A Riemannian metric defines a volume form, a Riemannian curvature and
a second fundamental form on $M$ (and, indeed, on each the $\partial_j
M$). Then
$\lips_{\dim M}(M)$ is its Riemannian volume and so is an integral of the
volume form. $\lips_{\dim M-1}(M)$ is closely related to the
(Riemannian) area of $\partial_{\dim M-1}(M)$, the ($\dim M
-1$)-dimensional boundary of $M$.
The remaining Lipschitz--Killing curvatures involve
integrals of both curvatures and second fundamental forms, and are more
complicated.
Rather than give the general form (for which you can see Chapter 10 of
\cite{TheBook}),
we shall give specific definitions for each of the cases we treat, in
the following
sections. Now, however, we shall give (relatively) simple expressions
for the cases
of interest to us. These are the main results of this paper.

\subsection{Main results}
\label{lkcrotatesubsec}

We shall present the results of this section in stages. In the first we shall
represent the Lipschitz--Killing curvatures as an integral for a
situation more general than
that needed to
handle rotation space fields. This is Theorem~\ref{decompositiontheorem},
and it shows how to write these integrals
as products of Lebesgue measures of subsets of $\real^N$ and integrals
over submanifolds of $\mathrm{Sym}_N^+$.
In Theorem~\ref{rotatetheorem} we shall compute these integrals for certain
situations, which will lead to Corollaries~\ref{rotatecorollary} and
\ref{corscale}. These give relatively simple, quite explicit, expressions
for the highest-order Lipschitz--Killing curvatures for rotation and
scale space fields, respectively.
More detailed information for scale space fields is given
in Section~\ref{scalesec}.

We say that a manifold
$L \subset\mathrm{Sym}_N^+$ is orthogonally invariant if
$S \in L \implies QSQ' \in L$ for every $Q \in O(N)$, the set of
orthogonal matrices.
For such $L$ of dimension~$m$, we define
%
\begin{equation}
\label{eqFdef} F_{m,j}(L) = \int_{L}
\mathcal{R}_j(S) \,d\mathcal{H}_m^{L}(S),
\end{equation}
where
$
\mathcal{R}_j(S)=\int_{O(N)}|(QSQ')_{j\times j}|^{
{1}/{2}}\,d\omega(Q)
$,
$\omega$ is the
Haar probability measure on $O (N )$, and for a matrix
$A$, $A_{j\times j}$ denotes the upper left submatrix of size $j\times j$.
The measure $\mathcal{H}_m^{L}$ is the Haar measure induced
by the restriction of the metric induced by the field $g$ to a fiber
${t} \times\mathrm{Sym}_N^+$, which, as we shall
see, is independent of the choice of $t$.
With these definitions, we can now state the following:
%
\begin{theorem}
\label{decompositiontheorem}
Let $f\dvtx A_l\times B_k\subset\mathbb{R}^N\times{\mathrm{Sym}}_N^+\to
\real$ be
the restriction
of the rotation space random field defined in Section \ref
{rotationspacesubsec}, to the product of an $l$-dimensional regular
stratified manifold $A_l\subset\mathbb{R}^N$ and a $k$-dimensional
regular stratified manifold $B_k\subset{\mathrm {Sym}}_N^+ $.
Furthermore, suppose that each stratum $\partial_jB_k$, $0\leq j\leq
k$, is orthogonally invariant and that the regularity conditions of
Section~\ref{regularitysubsec} vis-a-vis $f$ all hold. Then,
\begin{eqnarray*}
\lips_{l+k}(A_l \times B_k) &=&
\gamma^{l/2} \mathcal{H}_{l}^{E}(
\partial_l A_l) F_{k,l}(
\partial_k B_k),
\\
\lips_{l+k-1}(A_l
\times B_k)&=& \tfrac{1}{2} \bigl[\gamma^{(l-1)/2}
\mathcal{H}_{l-1}^E(\partial_{l-1}
A_l) F_{k,l-1}(\partial_k B_k)
\\
&&\hspace*{21.5pt}{} + \gamma^{l/2} \mathcal{H}_{l}^E(
\partial_l A_l) F_{k-1,l}(
\partial_{k-1} B_k) \bigr],
\end{eqnarray*}
where $\gamma$ is given by (\ref{ECscaledequn}) or (\ref{eqconstants}).
\end{theorem}

It is not hard to see that $\mathcal{R}_j(S)$ is a symmetric function
of the eigenvalues of $S$.
Furthermore, as we shall see, as each stratum of $B_m$ is
orthogonally invariant, its Hausdorff measure is in fact invariant
under $S\mapsto QSQ'$, for any orthogonal $Q$.
Hence, the $F_{m,j}$ are expressible as integrals over
certain functions of the eigenvalues of $S \in\partial_mB_m$.
In particular, for rotation space, with search region
$T \times\mathbb{S}(\sigma_2^{-2},\sigma_1^{-2})$,
we arrive at the following
explicit expressions.

\begin{theorem}
\label{rotatetheorem}
With the notation defined above, under the conditions of Theorem
\ref{decompositiontheorem}, and for a spherically symmetric kernel,
\begin{eqnarray*}
&&
F_{K,N}\bigl(\partial_K \mathbb{S}\bigl(
\sigma_2^{-2},\sigma_1^{-2}\bigr)
\bigr) \\
&&\qquad= C_{k,N} \int_{\sigma_{2}^{-2} \leq\lambda_1 \leq\cdots\leq\lambda
_{N} \leq\sigma_1^{-2}}\bigl\llvert \triangle
(\lambda )\bigr\rrvert \prod_{j=1}^{N}
\lambda_{j}^{-{N}/{2}}\,d\lambda,
\\
&&
F_{K,N-1}\bigl(\partial_K \mathbb{S}\bigl(
\sigma_2^{-2},\sigma_1^{-2}\bigr)
\bigr) \\
&&\qquad= C_{k,N} \int_{\sigma_{2}^{-2} \leq\lambda_1 \leq\cdots\leq\lambda
_{N} \leq\sigma_1^{-2}}\bigl\llvert \triangle
(\lambda )\bigr\rrvert \prod_{j=1}^{N}
\lambda_{j}^{-({N+1})/{2}}
\mathcal{R}_{N-1} \bigl(\operatorname{diag}(\lambda) \bigr)\,d
\lambda,
\\
&&
F_{K-1,N}\bigl(\partial_{K-1} \mathbb{S}\bigl(
\sigma_2^{-2},\sigma_1^{-2}\bigr)
\bigr) \\
&&\qquad= \frac{C_{k,N}}{D_{k,N}} \int_{\sigma_{2}^{-2} \leq\lambda_1 \leq\cdots\leq\lambda_{N-1}
\leq\sigma_1^{-2}}\bigl\llvert \triangle
(\lambda )\bigr\rrvert
\\
&&\hspace*{119pt}\qquad\quad{} \times \prod_{j=1}^{N-1}
\lambda_{j}^{-N/2}\Biggl[\sigma_{2}^{N-2}\prod
_{j=1}^{N-1}\bigl\llvert \lambda_{j}-
\sigma_{2}^{-2}\bigr\rrvert \\
&&\qquad\quad\hspace*{184.5pt}{} +\sigma_{1}^{N-2}
\prod_{j=1}^{N-1}\bigl\llvert
\sigma_{1}^{-2}-\lambda_{j}\bigr\rrvert \Biggr]\,d
\lambda,
\end{eqnarray*}
where $C_{k,N}$ and $D_{k,N}$ are given in (\ref{eqdefCk}) and (\ref
{eqdefDk}).
\end{theorem}
Theorem~\ref{decompositiontheorem} implies the
following corollary, with all $F_{m,j}$'s as above.
%
\begin{corollary}
\label{rotatecorollary}
Under the conditions of Theorem~\ref{rotatetheorem} and the
regularity conditions of Section
\ref{regularitysubsec},
the top two Lipschitz--Killing curvatures for the rotation random field
on $T \times\mathbb{S}(\sigma_2^{-2},\sigma_1^{-2})$ are
\begin{eqnarray*}
&&
\lips_{D}\bigl(T \times\mathbb{S}\bigl(
\sigma_2^{-2},\sigma_1^{-2}\bigr)
\bigr) \\
&&\qquad= 
\gamma^{N/2} \mathcal{H}_{N}^{E}(
\partial_NT) F_{K,N}\bigl(\partial_K\mathbb{S}
\bigl(\sigma_2^{-2},\sigma_1^{-2}
\bigr)\bigr),
\\
&&
\lips_{D-1}\bigl(T \times
\mathbb{S}\bigl(\sigma_2^{-2},\sigma_1^{-2}
\bigr)\bigr) 
\\
&&\qquad = \tfrac{1}{2} \bigl(\gamma^{(N-1)/2} \mathcal
{H}_{N-1}^E(\partial_{N-1} T) F_{K,N-1}
\bigl(\partial_K\mathbb{S}\bigl(\sigma_2^{-2},
\sigma_1^{-2}\bigr)\bigr)
\\
&&\hspace*{27pt}\qquad\quad{} + \gamma^{N/2} \mathcal{H}_{N}^E(
\partial_NT) F_{K-1,N}\bigl(\partial_{K-1}\mathbb{S}
\bigl(\sigma_2^{-2},\sigma_1^{-2}
\bigr)\bigr) \bigr).
\end{eqnarray*}
\end{corollary}

To obtain a parallel result for scale space fields, we
use the fact that the scale space random field on
$T\times[\sigma_1,\sigma_2]$, with a spherically symmetric kernel,
has the same distribution as the restriction of
the rotation space random field defined in
Section~\ref{rotationspacesubsec} to $T \times\mathbb{D}(\sigma_2^{-2},\sigma_1^{-2})$, where
\[
\mathbb{D}\bigl(\sigma_2^{-2}, \sigma_1^{-2}
\bigr) = \bigl\{ \nu I\dvtx  \sigma_2^{-2} \leq\nu\leq
\sigma_1^{-2} \bigr\} \subset\mathbb {S}\bigl(
\sigma_2^{-2}, \sigma_1^{-2}\bigr).
\]
Once we study the metric induced by the rotation field, Theorem \ref
{decompositiontheorem} will effortlessly yield the following result.
%
\begin{corollary}
\label{corscale}
Under the conditions of Theorem~\ref{rotatetheorem} and the
regularity conditions of Section
\ref{regularitysubsec},
the top two Lipschitz--Killing curvatures for the scale space random
field on
$T\times[\sigma_1,\sigma_2]$
are
%
\begin{eqnarray}
\label{scalelipsNequn} \lips_{N+1}\bigl(T
\times\mathbb{D}\bigl(\sigma_2^{-2},\sigma_1^{-2}
\bigr)\bigr) &=& \gamma^{N/2} \mathcal{H}_{N}^{E}(
\partial_NT) F_{1,N}\bigl(\partial_1\mathbb{D}
\bigl(\sigma_2^{-2},\sigma_1^{-2}
\bigr)\bigr)
\nonumber\\
&=& \widetilde{C}_{k,N} \gamma^{N/2} \mathcal{H}_{N}^{E}(
\partial_NT) \int_{\sigma_2^{-2}}^{\sigma_1^{-2}}
\nu^{N/2-1} \,d\nu
\\
&=& \widetilde{C}_{k,N} \gamma^{N/2} \mathcal{H}_{N}^{E}(
\partial_NT) \frac{1}{N/2} \bigl(\sigma_1^{-N}
- \sigma_2^{-N} \bigr)
\nonumber
\end{eqnarray}
and
%
\begin{eqnarray}\label{scalelipsNminus1equn}
&& \lips_{N}\bigl(T \times\mathbb{D}\bigl(
\sigma_2^{-2},\sigma_1^{-2}\bigr)
\bigr)
\nonumber\\
&&\qquad = \frac{1}{2} \bigl(\gamma^{(N-1)/2}\mathcal
{H}_{N-1}^E(\partial_{N-1} T) F_{1,N-1}
\bigl(\partial_1\mathbb{D}\bigl(\sigma_2^{-2},
\sigma_1^{-2}\bigr)\bigr)
\nonumber\\
&&\hspace*{49.5pt}\qquad\quad{} + \gamma^{N/2} \mathcal{H}_{N}^E(
\partial_NT) F_{0,N}\bigl(\partial_0
\mathbb{D}\bigl(\sigma_2^{-2},\sigma_1^{-2}
\bigr)\bigr) \bigr)
\\
&&\qquad = \frac{1}{2} \biggl(\widetilde{C}_{k,N}
\gamma^{(N-1)/2}\mathcal{H}_{N-1}^E(
\partial_{N-1} T) \frac{ (\sigma_1^{-(N-1)} - \sigma_2^{-(N-1)}  )}{(N-1)/2}
\nonumber
\\
&&\hspace*{101pt}\qquad\quad{} + \gamma^{N/2} \mathcal{H}_{N}^E(
\partial_NT) \bigl(\sigma_1^{-N} +
\sigma_2^{-N} \bigr) \biggr),
\nonumber
\end{eqnarray}
where $\widetilde{C}_{k,N}$ is as in (\ref{eqdefCtildek}).
\end{corollary}

In Section~\ref{scalesec} we analyze scale space fields independently
of the
rotation space case for a manifold with a boundary and without assuming
spherical symmetry. The formulae for all Lipschitz--Killing curvatures
are given via
(\ref{ECsumofLsequn}), the terms of which are all computed
in Section~\ref{scaleLKCSsubsubsec}.





\subsection{Comparison with existing results}

While the formulae of Theorem~\ref{rotatetheorem} may look rather
forbidding, one should
note two facts. The first is that, for dimension $N>2$, they are the
first of their
kind. The second is that although for dimension $N=2$ they are implicit in
\cite{Shafie03rotationspace}, one computation relied on computer
algebra and gave a large,
unstructured formula, while the second, based on tube formulae, does
not precisely give the
LKCs. (This is not unreasonable, since that computation was directed
explicitly at
approximating exceedance probabilities, and not computing LKCs.) To
recoup the
results of~\cite{Shafie03rotationspace}, we treat, as an example, the
case $N=2$ (so that $D=5$).
In that case, it is possible to
carry out the integrations needed to compute $\lips_5$ and $\lips_4$.
A page or so of calculus shows that, for the Gaussian kernel,
%
\begin{equation}
\lips_5(T \times\mathbb{S}) = 2^{-4}\pi
\sigma_{1}^{-2}\llvert T\rrvert \bigl[r^{2}-1-
\bigl(r^{2}+1 \bigr)\ln r \bigr],
\end{equation}
where $|T|$ is the simple Euclidean area of $T$ and $r=\sigma_{1}/\sigma_{2}$.

The LKC $\lips_4$ is a little harder to calculate, since while
the first term inside the brackets is again amenable to simple
calculus, the second term
leads us to special functions. However, a simple answer is attainable, and
\begin{eqnarray*}
\lips_{4} (T \times\mathbb{S} ) &=& 2^{-7/2}\pi
\sigma_{1}^{-2}\llvert T\rrvert \bigl(r^{2}-1
\bigr)\ln r
\\
&&{} +2^{-9/2}\pi\llvert \partial T\rrvert \int_{\sigma_{2}^{-2}\leq
\lambda_{1}\leq\lambda_{2}\leq\sigma_{1}^{-2}} (
\lambda_{2}-\lambda_{1} ) (\lambda_{1}
\lambda_{2} )^{-3/2}\mathcal{R} (\lambda_{1},
\lambda_{2} )\,d\lambda,
\end{eqnarray*}
where $|\partial T|$ is the length of $\partial T$,
and, with $E$ the elliptic integral
$E (y )=\int_{0}^{\pi/2}\sqrt{1-y^{2}\sin^{2}\theta
}\,d\theta$,
we have $\mathcal{R} (\lambda_{1},\lambda_{2} )={\pi
}\lambda_{1}
E ((\lambda_{1}^{2}-\lambda_{2}^{2})/\lambda_{1}^{2} )/2$.


Note also that, as stated,
Theorem~\ref{rotatetheorem} is weaker than the result for the scale
case in Section~\ref{scalesec}
in that it gives only the highest two LKCs, rather than all of them.
However, while more
would be desirable, it is the higher ones that are most important for
the Euler characteristic heuristic discussed in the following section.
We shall explain there why this is the case.


\section{The Euler characteristic heuristic: ECH}
\label{motivationsec}

As intimated in the \hyperref[secIntroduction]{Intro-} \hyperref[secIntroduction]{duction}, signal plus noise problems
often boil down to being able to compute the exceedance
probabilities (\ref{exceedence1equn}). For smooth Gaussian
processes, there is typically no way to
directly compute these probabilities, and they must be approximated, using
approximations that are accurate for the tail 5\% or so of the distribution.

There are a number of ways to carry out such approximations, and two of
these are at the center
of the papers~\cite{Shafie03rotationspace} and~\cite{S-WOR95},
discussed in the \hyperref[secIntroduction]{Introduction}. One relies on \textit{tube
formulae} and moves the problem from one about Gaussian processes on
Euclidean sets $M$
to one involving computation of volumes of tubes around embeddings of
$M$ in high-dimensional spheres.


The other is the Euler characteristic heuristic
\cite{Adler2000,TheBook} which provides a generally easier route,
along with the advantage that it also provides information on
sample paths at all levels. The ECH
states that for smooth random fields~$f$, on nice,
locally convex, compact sets $M$ (which will be the case in our
scenarios if $T$ is locally convex, since
$\mathbb D$ and $\mathbb S$ are convex)
and large $u$,
%
\begin{equation}
\label{approxequn} \P \Bigl\{\sup_{t\in M} f(t) \geq u \Bigr\}  \approx
\E \bigl\{\chi \bigl(A_u(f,M) \bigr) \bigr\}.
\end{equation}

In the case of smooth, Gaussian, $f$, with zero mean and constant unit
variance (but with no
stationarity assumptions) the approximation in (\ref{approxequn}) can
be quantified
(cf.~\cite{TheBook,TaylorTakemuraAdler2003Validity}), and it is
known that
%
\begin{equation}
\label{validitymainequation} \liminf_{u \rightarrow\infty}
-u^{-2} \log\bigl\llvert \operatorname {Diff}_{f,M}(u) \bigr
\rrvert  \geq \frac{1}{2} \biggl(1+ \frac{1}{ \sigma^2_c(f)} \biggr),
\end{equation}
where we write $ \operatorname{Diff}_{f,M}(u)$ for the difference of
the two terms in
(\ref{approxequn}), and
$\sigma^2_c(f)>0$ is an $f$-dependent constant, that is, in many cases,
computable.

The importance of (\ref{validitymainequation}) is that it shows that
the two terms in
(\ref{approxequn}), both of which are known from the general theory
of Gaussian processes
to be of order $u^{\dim M-1} \exp(-u^2/2)$, differ by a
term of exponentially smaller order, viz. of order $ \exp
(-u^2/2(1+\sigma^2_c(f))$.
This fact is what justifies the use of the ECH for Gaussian fields (in
which case it is no longer heuristic)
and in large part motivates the computations of this paper, and
also explains the claim of the previous section, that the leading
Lipschitz--Killing curvatures are
the important ones in applying the ECH.


\section{Before we start the proofs}
\label{beforesection}

As we wrote earlier, it should have been
possible to read the paper up to this point without a knowledge of
differential geometry.
This will not be true for the proofs.

The derivations of the leading Lipschitz--Killing curvatures for the
rotation space case, which is what we have concentrated on so far, and
all LKCs for the scale space case, which will be treated in Section
\ref{scalesec}, begin in the same fashion. Each starts by taking the
parameter space $M$, writing it as a stratified manifold, and carefully
identifying the strata $\partial_jM$ of (\ref{ppartitionequation}). The
next step will be to take the generic formulae for LKCs, given in
(necessary, but sometimes painful) detail in~\cite{TheBook}, and
rewrite them for the current setting, that is, for the specific
stratification, and for the induced Riemannian metric (\ref{metric}).
We shall not explain in detail how this is done other than to say, at
each stage, which generic formula we are using.
Then, the Riemannian
metric induced by the field is computed.

At this point the proofs diverge. For scale space
we compute all LKCs, and so geometric objects such as curvature and
second fundamental forms
arise. For rotation space, with only the two leading LKCs, we can avoid them.
The problem lies not in the the higher dimension
of rotation space, but in its more complicated Riemannian metric.
Nevertheless, in Section~\ref{secscale} we
shall see that once the geometry of rotation space is understood, the
two leading LKCs of scale space are effortlessly computed.

Throughout, we have tried to carry out computations as carefully as
possible, so that not only
can the interested reader follow the calculations, but also can see how
they need to be adjusted for
other cases.

%

\section{\texorpdfstring{Proof of Theorem \protect\ref{decompositiontheorem}}{Proof of Theorem 2.1}}
\label{rotatedecomposition}


We shall break the proof of Theorem~\ref{decompositiontheorem} into stages,
the first of which relates to the crucial
orthogonality properties mentioned in the \hyperref[secIntroduction]{Introduction}.

\subsection{\texorpdfstring{Orthogonality and the Riemannian metric on $\mathbb{R}^N\times\{S \}$}
{Orthogonality and the Riemannian metric on R N x \{S\}}}
\label{dependencesubsec}

Our aim in this subsection is to establish the
orthogonality of the two subspaces
\[
T_{t}\mathbb{R}^N, T_{S}
\mathrm{Sym}_{N}^+\subset T_{t}\mathbb {R}^N
\oplus T_{S}\mathrm{Sym}_{N}^+\simeq T_{ (t,S )}
\mathbb{R}^N\times{\mathrm{Sym}}_N^+,
\]
under the inner product defined by the Riemannian metric
induced by the process. Another
factor that comes into play is a structured dependence of the restrictions
of the Riemannian\vadjust{\goodbreak} metric on $\mathbb{R}^N\times{\mathrm{Sym}}_N^+$
to fibers of
the form $\mathbb{R}^N\times \{ S \} $.
This will turn out to be the key that allows
us to carry out further detailed computations.

To establish this independence, let $ (t,S )$ be a point in
$\mathbb{R}^N\times{\mathrm{Sym}}_N^+$ and let $X,Y$
be two derivations in $T_{t}\mathbb{R}^N, T_{S}{\mathrm
{Sym}}_N^+\subset
T_{ (t,S )}\mathbb{R}^N\times{\mathrm{Sym}}_N^+$,
respectively. Using the connection to the covariance function given
in (\ref{eqRiemannianmetric}), and writing
$Xt$ and $YS$ to denote the elementwise derivatives of the identity
functions on $\mathbb{R}^N$ and ${\mathrm{Sym}}_N^+$, respectively,
we have%
\begin{eqnarray*}
g (X,Y ) & = & X_{ (t,S )}Y_{ (\widetilde
{t},\widetilde{S} )}C\bigl((t,S),(\widetilde{t},\widetilde{S})
\bigr) \big|_{ (t,S )= (\widetilde{t},\widetilde{S} )}
\\[-2pt]
& = & 2\llvert S\rrvert^{{1}/{4}} \bigl(Y
\llvert S\rrvert^{{1}/{4}} \bigr) \int\bigl((Xt)'Su\bigr)
\dot{k}\bigl( u'Su\bigr) k\bigl(u'Su\bigr) \,du
\\[-2pt]
&&{} + 2\llvert S\rrvert^{{1}/{2}} 
\int \bigl((Xt)'Su
\bigr) \bigl(u'(YS)u\bigr) \bigl[\dot{k}\bigl( u'Su
\bigr)\bigr]^2 \,du\\[-2pt]
&=& 0,
\end{eqnarray*}
where the last equality follows since both integrands are odd functions.
This shows that the subspaces $T_{t}\mathbb{R}^N$ and $T_{S}{\mathrm
{Sym}}_N^+$
are orthogonal
at each point.

We now compute the restriction of the Riemannian metric $g$
to a fiber \mbox{$\mathbb{R}^N\times \{ S \} $}, again based on
(\ref
{eqRiemannianmetric}).
Writing, again, $\nabla_t f$ for
the spatial derivative of $f$ with respect to the
elements of $t$, the Gram matrix of the
Riemannian metric on \mbox{$\mathbb{R}^N\times \{ S \} $},
relative to the
standard Euclidean basis, is
\begin{eqnarray*}
&&\operatorname{Var} \bigl\{ \nabla_{t}f (t,S ) \bigr\}
\\[-2pt]
&&\qquad = \llvert S\rrvert^{{1}/{2}}\int_{\mathbb{R}^{N}}
\nabla_{t}k \bigl( (t-u )'S (t-u ) \bigr) \bigl(
\nabla_{t}k \bigl( (t-u )'S (t-u ) \bigr)
\bigr)'\,du
\\[-2pt]
&&\qquad = 4\llvert S\rrvert^{{1}/{2}}\int_{\mathbb{R}^{N}}\dot {k}
\bigl( w'Sw \bigr) Sw \bigl(\dot{k} \bigl( w'Sw \bigr)
Sw \bigr)'\,dw
\\[-2pt]
&&\qquad = 4\llvert S\rrvert^{{1}/{2}}S \biggl[\int_{\mathbb
{R}^{N}}
\bigl(\dot{k} \bigl( w'Sw \bigr) \bigr)^{2}ww'\,dw
\biggr]S
\\[-2pt]
&&\qquad = 4S^{{1}/{2}} \biggl[\int
_{\mathbb{R}^{N}} \bigl(\dot {k} \bigl( v'v \bigr)
\bigr)^{2}vv'\,dv \biggr]S^{{1}/{2}}
\\[-2pt]
&&\qquad = \gamma S.
\end{eqnarray*}


\subsection{\texorpdfstring{The Riemannian metric on $\{t\}\times\mathrm{Sym}^+_N$}
{The Riemannian metric on \{t\} x Sym+N}}

Relying on (\ref{eqRiemannianmetric}) to express the Riemannian
metric and using a simple change of variables, it can easily be seen
that the restriction of $g$ to $ \{ t \} \times\mathrm{Sym}_N^+$
is independent of $t$. Denote this metric on $\mathrm{Sym}_N^+\simeq
\{ t \} \times\mathrm{Sym}_N^+$
by $g^{\mathrm{Sym}_N^+}$.
Using the identity $\nabla_{S}\llvert S\rrvert =\llvert S\rrvert S^{-1}$,
valid on $\mathrm{Sym}_N^+$,
we have
\begin{eqnarray*}
&&{\nabla_{S} \bigl(\llvert S\rrvert^{{1}/{4}}k \bigl( (t-u
)'S (t-u ) \bigr) \bigr)}
\\[-2pt]
&&\qquad = \tfrac{1}{4}\llvert S\rrvert^{{1}/{4}}S^{-1}k \bigl(
(t-u )'S (t-u ) \bigr)
\\[-2pt]
&&\qquad\quad{} + \llvert S\rrvert^{{1}/{4}}\dot {k} \bigl( (t-u )'S (t-u
) \bigr) (t-u ) (t-u )'.
\end{eqnarray*}
Therefore, for tangent vectors $A, B \in\mathrm{Sym}_N \simeq
T_S\mathrm{Sym}^+_N$,
\begin{eqnarray*}
g_{S}^{\mathrm{Sym}_N^+} (A,B ) & = & \int_{\mathbb
{R}^{N}} \bigl
\langle\nabla_S \bigl(\llvert S\rrvert^{{1}/{4}}k \bigl( (t-u
)'S (t-u ) \bigr) \bigr),A \bigr\rangle
\\[-2pt]
&&\hspace*{15pt}{} \times\bigl\langle\nabla_S \bigl(\llvert S\rrvert^{{1}/{4}}k
\bigl( (t-u )'S (t-u ) \bigr) \bigr),B \bigr\rangle \,du
\\[-2pt]
&=& \llvert S\rrvert^{{1}/{2}}\int_{\mathbb{R}^{N}} \biggl[
\frac
{1}{4}k \bigl( u'Su \bigr) \bigl\langle
S^{-1},A \bigr\rangle+\dot {k} \bigl( u'Su \bigr) \bigl
\langle uu',A \bigr\rangle \biggr]
\\[-2pt]
&&\hspace*{43.4pt}{} \times\biggl[\frac{1}{4}k \bigl( u'Su \bigr) \bigl\langle
S^{-1},B \bigr\rangle+\dot{k} \bigl( u'Su \bigr) \bigl
\langle uu',B \bigr\rangle \biggr]\,du.
\end{eqnarray*}
%
After some slightly tedious calculus, and based on the fact that
$\kappa_4-3\kappa_{2,2}=0$ [see~(\ref{eqkappaidentity})], all of
the terms
above can be expressed in terms of
moments of the kernel $k$ and its derivative $\dot{k}$ as follows:
\begin{eqnarray*}
g_{S}^{\mathrm{Sym}_N^+} (A,B )
&=& \biggl(\frac{1}{16} + \frac{\rho}{2} + \kappa_{2,2}
\biggr) \Tr\bigl(S^{-1}A\bigr)\Tr\bigl(S^{-1}B\bigr)+ 2
\kappa_{2,2} \Tr \bigl(S^{-1}AS^{-1}B\bigr).
\end{eqnarray*}

When $h$ is the Gaussian kernel $h(u)=e^{-\|u\|^2_2/2}/\pi^{N/4}$, and
therefore
$k(x)=e^{-x/2}/\pi^{N/4}$, simple
calculations show that
$\frac{1}{16} + \frac{\rho}{2} + \kappa_{2,2}=0$, leaving only the term
$2\kappa_{2,2}\Tr(S^{-1}AS^{-1}B)$ above.

\subsection{Putting everything together}\label{puttogethersubsec}

Let $\mathcal{H}_{l}^{A_l \times\{S\}}$ denote the
Riemannian measures on $A_l \times \{ S \} \simeq A_l$
induced by $g$
and let $\mathcal{H}_{l}^{E}$ denote the Hausdorff measure, that is,
Riemannian measure under
the standard Euclidean metric on $T$.
A standard calculation gives %
%
\[
\frac{d\mathcal{H}_{l}^{A_l \times\{S\}}}{d\mathcal
{H}_{l}^{E}} (t )= \gamma^{l/2}\llvert S\rrvert_{T_{t}\partial_lA_l}^{{1}/{2}},
\]
where $\llvert S\rrvert_{T_{t}\partial_lA_l}$ is defined
to be the determinant of $PSP'$ for some $l\times N$
matrix $P$, the rows of which form an orthonormal basis %
of $T_{t}\partial_lA_l$ (under the Euclidean inner product).
The determinant, of course, is independent of the choice of basis.

Working in a chart and relying on the orthogonality established above,
it is
easy to see that
%
\begin{eqnarray}
\label{eq41} %
\lips_{k+l}
(A_l \times B_k ) &=& \hauss_{k+l}(A_l
\times B_k)
\nonumber\\[-2pt]
&=& \int
_{\partial_{k}B_k}\int_{\partial_{l}A_l}d\mathcal
{H}_{l}^{A_l \times\{S\}} (t )\,d\mathcal{H}_{k}^{B_k}
(S )
\\[-2pt]
&=& \gamma^{l/2} \int_{\partial_l A_l}\int
_{\partial_k B_k}\llvert S\rrvert_{T_{t}\partial_lA_l}^{{1}/{2}}\,d\mathcal
{H}_{k}^{B_k} (S )\,d\mathcal{H}_{l}^{E}
(t ). \nonumber
\end{eqnarray}
The equality above relies on the expression for the Radon--Nikodym
derivative above and
a Fubini theorem.

We recall now our assumption that $B_k$ is orthogonally invariant, that
is, $S \in B_k \implies QSQ' \in B_k$ for any $Q \in O(N)$. Thus, the\vadjust{\goodbreak}
mapping $S\mapsto QSQ'$ is a diffeomorphism from $B_k$
to itself. Moreover, for the covariance function of the rotation field,
we have
\[
C\bigl((t,S),(t,\widetilde{S})\bigr)=C\bigl(\bigl(t,QSQ'\bigr),
\bigl(t,Q\widetilde{S}Q'\bigr)\bigr).
\]
Therefore, $S\mapsto QSQ'$ is also a Riemannian isometry.
Thus, the measure $\hauss_k^{B_k}$, induced by the Riemannian metric,
is also
invariant under $S\mapsto QSQ'$.
This implies that we can drop the dependence of the integral in (\ref{eq41})
on the tangent space $T_{t}\partial_l A_l$. That is, we can replace
it, for all $t$, with a fixed $j$-dimensional subspace of our choice or
a random one.
Therefore,
%
\begin{eqnarray}
\hauss_{k+l}(A_l
\times B_k) 
&=& \gamma^{l/2} \hauss_l(\partial_lA_l)
\int_{\partial_{k}B_k} \mathcal{R}_l(S) \,d
\mathcal{H}_{k}^{B_k} (S )
\nonumber\\[-8pt]\\[-8pt]
&=& \gamma^{l/2} \hauss_l(\partial_lA_l)
F_{k,l}(\partial_kB_k). \nonumber
\end{eqnarray}

To complete the proof of the theorem, one must repeat
essentially the same calculation for the codimension
part of $A_l \times B_k$.

\section{\texorpdfstring{Proof of Corollary \protect\ref{corscale}}{Proof of Corollary 2.4}}
\label{secscale}

In this section we perform the computations necessary to
evaluate the expression in Corollary~\ref{corscale}.
These calculations are much simpler than those for rotation
space. The set $\mathbb{D}(\sigma_2^{-2}, \sigma_1^{-2})$
can be stratified as
\begin{eqnarray*}
\partial_1 \mathbb{D}\bigl(\sigma_2^{-2},
\sigma_1^{-2}\bigr)&=&\bigl\{\nu I\dvtx  \sigma_2^{-2}
< \nu< \sigma_1^{-2}\bigr\},\\
\partial_0\mathbb{D}\bigl(\sigma_2^{-2}, \sigma_1^{-2}
\bigr) &=& \bigl\{\sigma_2^{-2}I \bigr\} \cup\bigl\{
\sigma_1^{-2}I \bigr\}.
\end{eqnarray*}
Clearly, each stratum is orthogonally invariant. Consequently, Theorem
\ref{decompositiontheorem} implies we need only compute
the following three terms, each of which is straightforward:
%
%
%
\begin{eqnarray*}
F_{0,N}\bigl(\partial_0 \mathbb{D}\bigl(
\sigma_2^{-2}, \sigma_1^{-2}\bigr)
\bigr) &=& \sigma_1^{-N} + \sigma_2^{-N},
\\
F_{1,N}\bigl(
\partial_1 \mathbb{D}\bigl(\sigma_2^{-2},
\sigma_1^{-2}\bigr)\bigr)&=& \int_{\sigma_2^{-2}}^{\sigma_1^{-2}}
\nu^{N/2} f(\nu) \,d\nu,
\\
F_{1,N-1}\bigl(\partial_1 \mathbb{D}\bigl(\sigma_2^{-2},
\sigma_1^{-2}\bigr)\bigr)&=& \int_{\sigma_2^{-2}}^{\sigma_1^{-2}}
\nu^{(N-1)/2} f(\nu) \,d\nu, 
\end{eqnarray*}
where $f(\nu)$ is the appropriate density with respect to the Lebesgue measure.

Above, we expressed the necessary integrals in terms of the
parametrization of $\partial_1 \mathbb{D}(\sigma_2^{-2}, \sigma_1^{-2})$ by the curve $\nu\mapsto\nu\cdot I$.
Hence,
\begin{eqnarray*}
f(\nu) &=& \sqrt{g^{\mathrm{Sym}_N^+}_{\nu\cdot I}(I,I)}
\\
&=& \nu^{-1}\sqrt{N^2 \biggl(\frac{1}{16} +
\frac{\rho}{2} + \kappa_{2,2} \biggr) + 2N\kappa_{2,2}}
\\
& =& \widetilde{C}_{k,N} \nu^{-1}. 
\end{eqnarray*}
Therefore,
\begin{eqnarray*}
F_{1,N}\bigl(\partial_1 \mathbb{D}
\bigl(\sigma_2^{-2}, \sigma_1^{-2}
\bigr)\bigr)&=& \widetilde{C}_{k,N} \frac{1}{N/2} \bigl(
\sigma_1^{-N}-\sigma_2^{-N} \bigr),
\\
F_{1,N-1}\bigl(\partial_1 \mathbb{D}\bigl(\sigma_2^{-2},
\sigma_1^{-2}\bigr)\bigr)&=& \widetilde{C}_{k,N}
\frac{1}{(N-1)/2} \bigl(\sigma_1^{-(N-1)}-\sigma_2^{-(N-1)}
\bigr).
\end{eqnarray*}

\section{\texorpdfstring{Proof of Theorem \protect\ref{rotatetheorem}}{Proof of Theorem 2.2}}

We start the proof of Theorem~\ref{rotatetheorem} by describing a
stratification of $\mathbb{S}(\sigma_2^{-2},\sigma_1^{-2})$,
which is one of the regularity conditions required for applying the GKF.
In fact, we shall actually only use the strata of codimension 0 and 1.

\subsection{Stratification of the space}\label{secStratification}
Our first step is to express the manifold $M= T\times\mathbb
{S}(\sigma_2^{-2},\sigma_1^{-2})$
as a Whitney stratified manifold. Note that a stratification for
$\mathbb{S}(\sigma_2^{-2},\sigma_1^{-2})$ yields a stratification
for $M$ as
the product stratification, given that a Whitney stratification of $T$
is assumed
given. Thus, our first goal is to stratify
$\mathbb{S}(\tau_{1},\tau_{2})$
for arbitrary $0<\tau_{1}<\tau_{2}<\infty$.

Semialgebraic orbits of the Lie group action $GL_{N} (\mathbb
{R} )\times\mathrm{Sym}_{N}\rightarrow\mathrm{Sym}_{N}$,
$ (P,S )\mapsto PSP'$,
%
\begin{equation}\label{eqnonegdef-stratification}
\bigl\{ S\in\mathrm{Sym}_{N}| \operatorname{rank}(S )=k, S\geq0 \bigr\},\qquad
0\leq k\leq N,
\end{equation}
constitute a Whitney stratification of $ \{ S\in
\mathrm{Sym}_{N}| S\geq0 \} $
(cf.~\cite{Gibson}, page 21).

Since
\begin{eqnarray*}
A^{k} (\tau_{1} ) & \triangleq & \bigl\{ S\in
\mathrm{Sym}_{N} | \operatorname{rank} (S-\tau_{1}I )=N-k,
\tau_{1}I\leq S \bigr\},
\\
B^{k} (\tau_{2} ) & \triangleq & \bigl\{ S\in
\mathrm{Sym}_{N} | \operatorname{rank}(S-\tau_{2}I )=N-k, S\leq
\tau_{2}I \bigr\}
\end{eqnarray*}
are obtained by applying affine transformations to the strata of
(\ref{eqnonegdef-stratification}), they form Whitney stratifications
of
%
\begin{equation}\label{eq9}
\{ S\in\mathrm{Sym}_{N} | \tau_{1}I\leq S \} \quad\mbox{and}\quad
\{ S\in\mathrm{Sym}_{N} | S\leq\tau_{2}I \}.
\end{equation}

The final step in stratifying $\mathbb{S} (\tau_{1},\tau_{2} )$
is to verify that these stratifications intersect transversely. For
this, it suffices to show that for any $j,k\geq0$
such that $j+k\leq N$, a matrix $S\in A^{j}(\tau_{1})\cap B^{k}(\tau_{2})$
and nonzero matrix $C\in\mathrm{Sym}_{N}$, there exist smooth paths
$A (t )\dvtx  (-\delta,\delta )\longrightarrow
A^{j}(\tau_{1})$,
$B (t )\dvtx  (-\delta,\delta )\longrightarrow
B^{k}(\tau_{2})$
with $A (0 )=B (0 )=S$ such that for any $t\in
(-\delta,\delta )$,
%
\begin{equation}\label{eq3}
A (t )-B (t )= \bigl(A (t )-S \bigr)- \bigl(B (t )-S \bigr)=tC.
\end{equation}

Based on the fact that the group action above preserves
rank, the problem of finding the paths $A (t ), B
(t )$
is equivalent to defining them by
%
\begin{eqnarray}
\label{eq4}
A (t )&=& X (t ) (S-\tau_{1}I ) \bigl(X (t ) \bigr)'+
\tau_{1}I,
\\
\label{eq5}
B (t )&=& Y (t ) (S-\tau_{2}I ) \bigl(Y (t ) \bigr)'+
\tau_{2}I
\end{eqnarray}
and finding appropriate invertible matrices $X (t ),
Y (t )$.\vadjust{\goodbreak}

Even under the restriction $X (t )=Y (t )$ a solution
exists. Substitution in (\ref{eq3}) gives
\[
X (t ) (S-\tau_{1}I ) \bigl(X (t ) \bigr)'+
\tau_{1}I-X (t ) (S-\tau_{2}I ) \bigl(X (t )
\bigr)'-\tau_{2}I=tC,
\]
which is equivalent to
\[
X (t ) \bigl(X (t ) \bigr)'=\frac{t}{\tau_{2}-\tau_{1}}C+I.
\]

For small $t$, define $X (t )$ to be the unique positive
definite square root of $\frac{t}{\tau_{2}-\tau_{1}}C+I$. Substituting
$X (t )$ in (\ref{eq4}) and (\ref{eq5}) defines, for
small $\delta$, the required paths $A (t )$, $B
(t )$,
and proves transversality.

The transversality of the intersections, together with the stratifications
in (\ref{eq9}), yield the Whitney stratification
\[
\mathbb{S} (\tau_{1},\tau_{2} )=\mathop{\bigcup
_{0\leq j,k}}_{j+k\leq N}\mathbb{S}^{j,k} (
\tau_{1},\tau_{2} ),
\]
where $\mathbb{S}^{j,k} (\tau_{1},\tau_{2} )=A^{j}
(\tau_{1} )\cap B^{k} (\tau_{2} )$.
That is, we have found that the partition of $\mathbb{S} (\tau_{1},\tau_{2} )$
according to the geometrical multiplicity of $\tau_{1}$ and $\tau_{2}$
is a Whitney stratification.
%

Note for later use that the submanifold of $\mathrm{Sym}_{N}$ of
matrices of rank $N-k$, and
thus $A^k(\tau_{1})$ and $B^k(\tau_{2})$, are of codimension
$ k (k+1 )/2$. Since, for $j,k\geq0$ with $j+k\leq N$,
$A^j(\tau_{1})$ and $B^k(\tau_{2})$ intersect tranversely,
the codimension of $\mathbb{S}^{j,k}(\tau_{1},\tau_{2})$ is
$[j (j+1 )+ k (k+1 )]/2$.

\subsection{Volume density on $\mathrm{Sym}_N^+$}

Having computed the Riemannian metric on $\mathrm{Sym}_N^+$, in order
to apply Theorem~\ref{decompositiontheorem} to the rotation space
random field on $T \times\mathbb{S}(\sigma_2^{-2},\sigma_1^{-2})$,
we first\vspace*{1pt} explicitly compute $d\hauss_K^{\mathbb{S}(\sigma
_2^{-2},\sigma_1^{-2})}/d\hauss_K^E$, that is, the density with
respect to the Lebesgue measure on\vspace*{1pt}
$\mathrm{Sym}_N^+$ when viewed as $\real^K$ via
the choice of a basis for the Hilbert--Schmidt inner product on $\mathrm{Sym}_N$.

Before turning to the calculation we
need some notation. Denote the standard basis of $\mathrm{Sym}_{N}$
by $e= (e_{ij} )_{i=1,j=i}^{N,N}$. Normalizing the elements
$e_{ij}$ with $i\neq j$ by a factor of ${1}/{\sqrt{2}}$ yields
an orthonormal basis relative to the Hilbert--Schmidt inner product
which we denote by $\widetilde{e}= (\widetilde{e}_{ij}
)_{i=1,j=i}^{N,N}$.

Working with the chart defined by mapping a matrix to its coordinates
relative to $\widetilde{e}$, we see that the density is given by
%
\begin{equation}
\label{eqdensity}\frac{d\mathcal{H}_{K}^{\mathbb{S} (\sigma
_{2}^{-2},\sigma_{1}^{-2} )}}{d\mathcal{H}_{K}^{E}} (S
)=\sqrt{\operatorname{det}\bigl\{ g_{S}^{\mathrm{Sym}_{N}^{+}}
(\widetilde{e}_{m}, \widetilde{e}_{n} ) \bigr\}_{m=1,n=1}^{K,K}},
\end{equation}
where $ (\widetilde{e}_{m} )_{m=1}^{K} =  (\widetilde
{e}_{i(m),j(m)} )_{m=1}^{K}$ is obtained
by vectorizing $ (\widetilde{e}_{ij} )_{i=1,j=i}^{N,N}$.

As already stated in Section~\ref{puttogethersubsec}, for any
$Q\in O (N )$,
the mapping $S\mapsto QSQ'$ is a Riemannian isometry of $g_{S}^{\mathrm{Sym}_{N}^{+}}$,
and thus preserves the measure $\mathcal{H}_{K}^{\mathbb{S}
(\sigma_{2}^{-2},\sigma_{1}^{-2} )}$.
Clearly, it also preserves the Euclidean Hausdorff measure $\mathcal
{H}_{K}^{E}$. Therefore,
the density in (\ref{eqdensity}) is invariant under $S\mapsto QSQ'$
as well. Since any real symmetric
matrix is orthogonally diagonalizable, it follows that the density only
depends on the eigenvalues of $S$. It is therefore sufficient to
compute it for diagonal matrices.

Let $\Lambda= \operatorname{diag}(\lambda_1,\ldots,\lambda_N
)$. For any two basis elements of
$\widetilde{e}$,
\begin{eqnarray*}
g_{\Lambda}^{\mathrm{Sym}_N^+} (\widetilde{e}_{ij},\widetilde
{e}_{kl} ) &=& \bigl(\tfrac{1}{16}+\kappa_{2,2}+
\tfrac{1}{2}\rho \bigr)\operatorname{Tr} \bigl(\Lambda^{-1}
\widetilde{e}_{ij} \bigr)\operatorname{Tr} \bigl(\Lambda^{-1}
\widetilde{e}_{kl} \bigr)
\\
&&{} +2\kappa_{2,2}\operatorname{Tr} \bigl(\Lambda^{-1}
\widetilde{e}_{ij}\Lambda^{-1}\widetilde{e}_{kl}
\bigr)
\\
&=& \bigl(\tfrac{1}{16}+\kappa_{2,2}+\tfrac{1}{2}\rho
\bigr)\lambda_{i}^{-1}\lambda_{k}^{-1}
\delta_{ij}\delta_{kl} +2\kappa_{2,2}
\lambda_{i}^{-1}\lambda_{j}^{-1}
\delta_{ik}\delta_{jl}. %
\end{eqnarray*}

Suppose we order the basis in such a way that its first $N$ elements
are $\widetilde{e}_{1,1},\ldots,\widetilde{e}_{N,N}$.
Then, the corresponding $K \times K$ matrix has the form
\begin{eqnarray*}
\bigl\{ g_{\Lambda}^{\mathrm{Sym}_N^+} (\widetilde{e}_{m},
\widetilde{e}_{n} ) \bigr\}_{m=1,n=1}^{K,K} &=& 2
\kappa_{2,2} \operatorname{diag}\bigl( \bigl\{ \lambda_{i (m
)}^{-1}
\lambda_{j (m )}^{-1} \bigr\}_{m=1}^{K}
\bigr)
\\
&&{} + \bigl(\tfrac{1}{16}+\kappa_{2,2}+\tfrac{1}{2}\rho
\bigr)VV',
\end{eqnarray*}
where
$V=(\{ \lambda_{i}^{-1}\}_{i=1}^{N},0,\ldots,0)'\in\real^K$.

By inspection, this is a rank-one perturbation of a diagonal matrix.
Now, for a diagonal matrix $D$ and vector $u$, $\operatorname{det}
(D+uu' )=\operatorname{det} (D ) (1+u'D^{-1}u )$,
which, when applied in this setting, yields
\[
\frac{d\hauss_K^{\mathbb{S}(\sigma_2^{-2},\sigma_1^{-2})}}{d\hauss_K^E}(\Lambda)
= \sqrt{\operatorname{det} \bigl\{ g_{\Lambda}^{\mathrm{Sym}_N^+}
(\widetilde{e}_{m}, \widetilde{e}_{n} ) \bigr\}_{m=1,n=1}^{K,K}} =
c_{k,N} \llvert \Lambda\rrvert^{{- (N+1 )}/{2}},
\]
where $c_{k,N}$ is defined in (\ref{eqdefck}). By the argument above,
for general $S \in\mathbb{S}(\sigma_2^{-2},\allowbreak\sigma_1^{-2})$,
\[
\frac{d\hauss_K^{\mathbb{S}(\sigma_2^{-2},\sigma_1^{-2})}}{d\hauss_K^E}(S)
= c_{k,N} \llvert S\rrvert^{{- (N+1 )}/{2}}.
\]

\subsection{The push-forward to eigenvalues}

Rewriting the integral defining $F_{K,N}(\mathbb{S}(\sigma_2^{-2},\sigma_1^{-2}))$,
(\ref{eqFdef}), in terms of the density
${d\mathcal{H}_{K}^{\mathbb{S}(\sigma_2^{-2},\sigma
_1^{-2})}}/{d\hauss_K^E}$, gives a relatively complicated
integral. Fortunately, there is a simple way to convert it to an integral
against a simpler density and, more importantly, over a much simpler
domain. To do this, we need first to recall a classical result from
the theory of random matrices.

The starting point of the computation of the
joint density of the eigenvalues
in the
Gaussian orthogonal ensemble (GOE) model (cf.~\cite{RM}) is similar to
our situation. We begin with a measure $\mu$ on $\mathrm{Sym}_{N}$ with\vadjust{\goodbreak}
density (relative to $\hauss_K^E$) that depends strictly on the
eigenvalues. The goal is to compute the density of the push-forward of
$\mu$, denoted by $\nu$, under the mapping $\alpha\dvtx
\mathrm{Sym}_{N}\rightarrow\mathbb{R}^{N}$ which maps a matrix to the
vector composed of its ordered eigenvalues.

In the case of the GOE, the density of $\mu$ is simply the product
of the densities of the independent entries of the (real symmetric)
random matrix. The entries of the matrix are zero-mean Gaussian variables
with variance $2$ on the diagonal and variance $1$ off the diagonal.
That is,
\begin{eqnarray*}
\frac{d\mu}{d\hauss_K^E} (S )&=& 2^{N-{D}/{2}}2^{-{N}/{2}} (2\pi
)^{-{N
(N+1 )}/{4}}\exp \Biggl\{ -\frac{1}{4}\sum
_{i=1}^{N}S_{ii}-\frac{1}{2}\sum
_{1\leq i<j\leq N}S_{ij} \Biggr\}
\\
&=& 2^{({N-D})/{2}} (2\pi )^{-{N (N+1
)}/{4}}\exp \biggl\{ -\frac{1}{4}
\operatorname{Tr} \bigl(S^{2} \bigr) \biggr\},
\end{eqnarray*}
where the factor of $2^{N-{D}/{2}}$ comes from the fact that
$\hauss_K^E$ was defined via identification
with the Hilbert--Schmidt basis and not the standard basis of symmetric matrices.

The classical result is that
\[
\frac{d\nu}{dl} (\lambda )= (2\pi )^{-
{N}/{2}}2^{-{N (N+1 )}/{4}}\prod
_{j=1}^{N}\frac{\Gamma
(1/2 )}{\Gamma (j/2 )}
\mathbf{1}_{\lambda
_{1}\leq\cdots\leq\lambda_{N}}\bigl\llvert \triangle (\lambda )\bigr\rrvert \exp
\Biggl\{ -\frac{1}{4}\sum_{i=1}^{N}
\lambda_{i}^{2} \Biggr\},
\]
where $l$ is the Lebesgue measure on $\mathbb{R}^{N}$ and $\triangle
(\lambda )\triangleq\prod_{i<j} (\lambda_{j}-\lambda_{i} )$
is the Vandermonde determinant.

The specific Gaussian nature of $\mu$ has only a limited role to
play in this result, and all that is relevant for our needs is that
the density is determined
by $\alpha (S )$ at each point. In fact, for any measure
$\widetilde{\mu}$ with such a density, we can write ${d\widetilde
{\mu}}/{d\hauss_K^E} (S )=\zeta_{\widetilde{\mu}}
(\alpha (S ) )$
and rely on the GOE argument to conclude that, for the corresponding
push-forward measure $\widetilde{\nu}$,
\begin{eqnarray*}
\frac{d\widetilde{\nu}}{dl} (\lambda ) &=& \frac{d\nu}{dl} (\lambda ) \biggl(
\frac{\zeta_{\widetilde{\mu}}}{\zeta_{\mu}} \biggr) (\lambda ) \\
&=& (2\pi)^{{D}/{2}-N} \prod
_{j=1}^{N}\frac{\Gamma
(1/2 )}{\Gamma (j/2 )}\mathbf{1}_{\lambda_{1}\leq
\cdots\leq\lambda_{N}}
\bigl\llvert \triangle (\lambda )\bigr\rrvert \zeta_{\widetilde{\mu}} (\lambda ),
\end{eqnarray*}
where $\mu$ and $\nu$ are the measures corresponding to the GOE model.

In particular, denoting the push-forward of $\mathcal{H}_{K}^{\mathbb
{S}(\sigma_2^{-2},\sigma_1^{-2})}$
by $\alpha$ by $\rho$ gives
\begin{eqnarray*}
\frac{d\rho}{dl} (\lambda ) & = & (2\pi )^{{D}/{2}-N}c_{k,N}\prod
_{j=1}^{N} \biggl(\frac{\Gamma
(1/2 )}{\Gamma (j/2 )}
\lambda_{j}^{-
({N+1})/{2}} \biggr)\mathbf{1}_{\lambda_{1}\leq\cdots\leq\lambda
_{N}}\bigl
\llvert \triangle (\lambda )\bigr\rrvert
\\
& = & C_{k,N}\mathbf{1}_{\lambda_{1}\leq\cdots\leq\lambda_{N}}\bigl\llvert \triangle (\lambda )
\bigr\rrvert \prod_{j=1}^{N}
\lambda_{j}^{-({N+1})/{2}}.
\end{eqnarray*}

\subsection{The two leading LKCs}

In view of Theorem~\ref{decompositiontheorem}, having computed the
density $\frac{d\rho}{dl} (\lambda )$,
the computation of $\mathcal{L}_{D} (M )$ simply becomes
a matter of rewriting the integral of $F_{K,N} (\partial_{K}\mathbb{S} (\sigma_{2}^{-2},\sigma_{1}^{-2} )
)$ in terms of
eigenvalues. Doing so leads to
\begin{eqnarray*}
F_{K,N} \bigl(\partial_{K}\mathbb{S} \bigl(
\sigma_{2}^{-2},\sigma_{1}^{-2} \bigr)
\bigr)& = &\int_{\partial_{K}\mathbb{S}
(\sigma_{2}^{-2},\sigma_{1}^{-2} )}\llvert S\rrvert^{
{1}/{2}}\,d
\mathcal{H}_{K}^{\partial_{K}\mathbb{S} (\sigma
_{2}^{-2},\sigma_{1}^{-2} )}
\\
& =& C_{k,N}\int_{\sigma_{2}^{-2}\leq\lambda_{1}\leq\cdots\leq
\lambda_{N}\leq\sigma_{1}^{-2}}\bigl\llvert \triangle (
\lambda )\bigr\rrvert \prod_{j=1}^{N}
\lambda_{j}^{-{N}/{2}}\,d\lambda.
\end{eqnarray*}

We now turn to $F_{K,N-1}(\partial_K\mathbb{S}(\sigma_2^{-2},\sigma_1^{-2}))$,
for which an identical argument yields
\begin{eqnarray*}
&&F_{K,N-1}\bigl(\partial_K\mathbb{S}\bigl(
\sigma_2^{-2},\sigma_1^{-2}\bigr)
\bigr)
\\
&&\qquad = C_{k,N} \int_{{\sigma_{2}^{-2} \leq\lambda_1 \leq
\cdots\leq\lambda_{N} \leq\sigma_1^{-2}}}\bigl\llvert \triangle (
\lambda )\bigr\rrvert \prod_{j=1}^{N}
\lambda_{j}^{-
({N+1})/{2}}\mathcal{R}_{N-1}\bigl(
\operatorname{diag}(\lambda)\bigr) \,d\lambda.
\end{eqnarray*}


It remains to evaluate $F_{K-1,N}(\mathbb{S}(\sigma_2^{-2},\sigma_1^{-2}))$. Although
this depends on $\mathcal{H}_{K-1}^{\partial_{K-1} \mathbb{S}
(\sigma_{2}^{-2},\sigma_{1}^{-2} )}$, we shall not have
to study this measure directly in order to compute the integral.
Rather, an application
of Federer's coarea formula \cite{TheBook,Federer} with a carefully
chosen map will suffice.

For any $0<\tau<\sigma_{1}^{-2}$ define
$\mathbb{S}_{\mathrm{min}}^u(\tau,\sigma_{1}^{-2})$ to be the set of real,
symmetric matrices with unique minimal eigenvalue and with all
eigenvalues in $ (\tau,\sigma_{1}^{-2} )$.
Consider the smooth %
mapping $\lambda_{\mathrm{min}}\dvtx \mathbb{S}_{\mathrm{min}}^u(0,\sigma_{1}^{-2})\rightarrow\mathbb{R}$
which maps a matrix to its minimal eigenvalue. To apply the
coarea formula, we need a Riemannian structure on the domain and range
of $\lambda_{\mathrm{min}}$. For the domain we take $g_{S}^{\mathrm{Sym}_{N}^{+}}$
as the Riemannian metric, and on $\mathbb{R}$ consider the Euclidean
metric.

Now, fix some $0<\tau<\sigma_{1}^{-2}$. Let $S\in\mathbb
{S}_{\mathrm{min}}^u(\tau,\sigma_{1}^{-2})$
and write $S=Q\Lambda Q'$ with $Q\in O (N )$,
$\Lambda=\operatorname{diag} (\lambda )$
diagonal and $\lambda=(\lambda_1,\ldots,\lambda_N)$ ordered so that
$\lambda_{1}\leq\cdots\leq\lambda_{N}$.

Defining $\widetilde{e}{}^{Q}_{ij}=Q\widetilde{e}_{ij}Q'$, $\widetilde
{e}{}^{Q}= (\widetilde{e}{}^{Q}_{ij} )_{i=1,j=i}^{N,N}$ forms a
basis of
$\mathrm{Sym}_{N}^{+}$.
%
%
Let $\mathrm{Proj}_V$ be the orthogonal projection onto
\[
V \triangleq\operatorname{span} \bigl( \bigl\{ \widetilde{e}_{ij}^{Q}
\bigr\}_{i=1,j=i}^{N,N}\setminus\bigl\{ \widetilde{e}_{11}^{Q}
\bigr\} \bigr).
\]
Let $ (b_{ij} )_{i=1,j=i}^{N,N}$ be an orthonormal
basis of $T_{S}\mathbb{S}_{\mathrm{min}}^u (\tau,\sigma_{1}^{-2} )$
relative to $g_{S}^{\mathrm{Sym}_{N}^{+}}$, such that for all $b_{ij}$
with $ (i,j )\neq(1,1 )$, $b_{ij} \in V$, and such
that, denoting $\bar{b}_{11} \triangleq
\widetilde{e}_{11}^{Q}-\mathrm{Proj}_V ( \widetilde{e}_{11}^{Q}
)$,
\[
b_{11}= \bar{b}_{11}/\Vert\bar{b}_{11} \Vert
= (\lambda_1/D_{k,N}) \bar{b}_{11},
\]
with $D_{k,N}$ given in (\ref{eqdefDk}). The last equality follows
from a page or two of calculus, using the projection formula of
\cite{TheBook}, page 173, and the Sherman-Morrison formula (see
\cite{Press2007}, Section 2.7.1).\vadjust{\goodbreak}

Since the basis is orthonormal, the Jacobian can be expressed by
\[
J\lambda_{\mathrm{min}} (S )= \Biggl[\sum_{i=1,j=i}^{N,N}
\bigl\langle\nabla\lambda_{\mathrm{min}} (S ),b_{ij} \bigr
\rangle^{2} \Biggr]^{1/2}.
\]

Note that if we can show that $ \langle\nabla\lambda_{\mathrm{min}}
(S ),\widetilde{e}_{ij}^{Q} \rangle$
is zero for any $ (i,j )\neq(1,1 )$, then it will
follow that
\[
J\lambda_{\mathrm{min}} (S )= (\lambda_1/D_{k,N} ) \bigl
\vert\bigl\langle\nabla\lambda_{\mathrm{min}} (S ),\widetilde
{e}_{11}^{Q} \bigr\rangle\bigr\vert.
\]
In
order to show this, and also compute this value, we
examine the eigenvalues of $S+t\widetilde{e}_{ij}^{Q}$,
for small $t$. Write $ S+t\widetilde{e}_{ij}^{Q}=Q
(\Lambda+t\widetilde{e}_{ij} )Q', $ and note that therefore we can work
with the eigenvalues of $U_{t}\triangleq\Lambda+t\widetilde{e}_{ij}$
instead. When $i=j$, $U_{t}$ is diagonal with eigenvalues $ \{
\lambda_{1},\ldots,\lambda_{i-1},\lambda_{i}+t,\lambda_{i+1},\ldots,\lambda_{N}
\} $.
When $i\neq j$, $U_{t}$ also takes a very simple form, with 
eigenvalues
%
\begin{eqnarray*}
&& \{ \lambda_{1},\ldots,\lambda_{i-1},
\lambda_{i+1},\ldots,\lambda_{j-1},\lambda_{j+1},\ldots,\lambda_{N} \}
\\
&&\qquad{} \cup \bigl\{ \tfrac{1}{2} \bigl(\lambda_{i}+
\lambda_{j}\pm\sqrt{ (\lambda_{i}-\lambda_{j}
)^{2}+2t^{2}} \bigr) \bigr\}.
\end{eqnarray*}
Therefore,
\[
\bigl\langle\nabla\lambda_{\mathrm{min}} (S ),\widetilde {e}_{ij}^{Q}
\bigr\rangle=\lim_{t\rightarrow0}\frac{\lambda_{\mathrm{min}}
(S+t\widetilde{e}_{ij}^{Q} )-\lambda_{\mathrm{min}} (S )}{t}=\delta_{
(1,1 )} (i,j ),
\]
which implies
$
J\lambda_{\mathrm{min}} (S )=\frac{\lambda_{\mathrm{min}} (S )}{D_{k,N}}
$.
Finally, Federer's coarea formula gives
\begin{eqnarray*}
&& \int_{\mathbb{S}_{\mathrm{min}}^{u} (\tau,\sigma_{1}^{-2}
)}J\lambda_{\mathrm{min}} (S )\llvert S
\rrvert^{
{1}/{2}}\,d\mathcal{H}_{K}^{\mathbb{S}_{\mathrm{min}}^{u} (\tau,\sigma
_{1}^{-2} )} (S )
\\
&&\qquad = \int_{ (\tau,\sigma
_{1}^{-2} )}\,dx\int_{\lambda_{\mathrm{min}}^{-1} (x )}\llvert S
\rrvert^{{1}/{2}}\,d\mathcal{H}_{K-1}^{\mathbb
{S}_{\mathrm{min}}^{u} (\tau,\sigma_{1}^{-2} )} (S ).
\end{eqnarray*}
Since the integral on $\lambda_{\mathrm{min}}^{-1} (x )$ is continuous
in $x$, differentiating by $\tau$ yields
\begin{eqnarray*}
&& \int_{\lambda_{\mathrm{min}}^{-1} (\sigma_{2}^{-2} )}\llvert S\rrvert^{{1}/{2}}\,d
\mathcal{H}_{K-1}^{\mathbb
{S}_{\mathrm{min}}^{u} (\tau,\sigma_{1}^{-2} )} (S )
\\
&&\qquad =-\frac{d}{d\tau}\bigg|_{\tau=\sigma _{2}^{-2}}\int_{\mathbb{S}_{\mathrm{min}}^{u}
(\tau,\sigma _{1}^{-2} )}J \lambda_{\mathrm{min}} (S )\llvert
S\rrvert^{{1}/{2}}\,d\mathcal{H}_{K}^{\mathbb{S}_{\mathrm{min}}^{u}
(\tau,\sigma_{1}^{-2} )} (S )
\\
&&\qquad =-\frac{d}{d\tau}\bigg|_{\tau=\sigma
_{2}^{-2}}\frac{C_{k,N}}{D_{k,N}}\int
_{\sigma_{2}^{-2}\leq\lambda
_{1}\leq\cdots\leq\lambda_{N}\leq\sigma_{1}^{-2}}\bigl\llvert \triangle (\lambda )\bigr\rrvert
\lambda_{1}\prod_{j=1}^{N}
\lambda_{j}^{-{N}/{2}}\,d\lambda
\\
&&\qquad =\frac{C_{k,N}}{D_{k,N}}\sigma_{2}^{N-2}\int
_{\sigma
_{2}^{-2}\leq\lambda_{1}\leq\cdots\leq\lambda_{N-1}\leq\sigma
_{1}^{-2}}\bigl\llvert \triangle (\lambda )\bigr\rrvert \prod
_{j=1}^{N-1} \bigl(\lambda_{j}^{-{N}/{2}}
\bigl\llvert \lambda_{j}-\sigma_{2}^{-2}\bigr
\rrvert \bigr)\,d\lambda,
\end{eqnarray*}
where the second equality is implied by the fact that any point of
\[
\bigl\{ \lambda\dvtx \tau\leq\lambda_{1}\leq\cdots\leq
\lambda_{N}\leq\sigma_{1}^{-2} \bigr\} \setminus
\alpha \bigl(\mathbb {S}_{\mathrm{min}}^{u} \bigl(\tau,
\sigma_{1}^{-2} \bigr) \bigr)
\]
is a boundary point or has zero density and the last equality follows
by writing the integral as an iterated integral and differentiating.

Obviously, we can also consider the maximal eigenvalue, define
$\lambda_{\mathrm{max}}$ and $\mathbb{S}_{\mathrm{max}}^{u} (\sigma_{2}^{-2},\tau )$
accordingly, and carry out a similar computation. Combined with the
last result and the fact that $\partial_{K-1}\mathbb{S}
(\sigma_{2}^{-2},\sigma_{1}^{-2} )=\lambda_{\mathrm{min}}^{-1}
(\sigma_{2}^{-2} )\cup\lambda_{\mathrm{max}}^{-1} (\sigma_{1}^{-2} )$, we can
compute the final term needed to complete the proof:
\begin{eqnarray*}
&& F_{K-1,N} \bigl(\partial_{K-1} \mathbb{S} \bigl(
\sigma_{2}^{-2},\sigma_{1}^{-2} \bigr)
\bigr)
\\[-2pt]
&&\qquad =\int_{\partial_{K-1}\mathbb{S} (\sigma
_{2}^{-2},\sigma_{1}^{-2} )}\llvert S\rrvert^{
{1}/{2}}\,d
\mathcal{H}_{K-1}^{\mathbb{S} (\sigma_{2}^{-2},\sigma
_{1}^{-2} )} (S )
\\[-2pt]
&&\qquad =\frac{C_{k,N}}{D_{k,N}}\int_{\sigma_{2}^{-2}\leq\lambda
_{1}\leq\cdots\leq\lambda_{N-1}\leq\sigma_{1}^{-2}}\bigl\llvert \triangle (
\lambda )\bigr\rrvert
\\[-2pt]
&&\hspace*{120.4pt}\qquad\quad{} \times \prod_{j=1}^{N-1}
\lambda_{j}^{-{N}/{2}}\Biggl[\sigma_{2}^{N-2}\prod
_{j=1}^{N-1}\bigl\llvert \lambda_{j}-
\sigma_{2}^{-2}\bigr\rrvert\\[-2pt]
&&\qquad\quad\hspace*{186pt}{} +\sigma_{1}^{N-2}
\prod_{j=1}^{N-1}\bigl\llvert
\sigma_{1}^{-2}-\lambda_{j}\bigr\rrvert \Biggr]\,d
\lambda.
\end{eqnarray*}
%

\section{Scale space fields}
\label{scalesec}

In this section we shall compute \textit{all} the Lipschitz--Killing
curvatures for the scale
space random
field, under the assumptions of Section~\ref{scalespacesubsec} on the
kernel $h$ and of Section
\ref{regularitysubsec} on the parameter space $T$.
Also, we require only that $h$ satisfy
(\ref{ECscaledequn}), and not full rotational symmetry.

In the final analysis, the formulae that we shall obtain are quite
simple. They involve no more than the Euclidean Lipschitz--Killing
curvatures of $T$, the parameters $\gamma$ and $\kappa$ of
(\ref{ECscaledequn}) and (\ref{ECscale-variancesequn}), and
combinatorial coefficients. On the other hand, there are many elements
to each formula, so that any explicit computation of them will involve
a computer. Thus, we shall proceed by setting out two general
formulae---(\ref {eqdefcurvaturemeasureonetubes}) and
(\ref{ECsumofLsequn})---and then in Section~\ref{scaleLKCSsubsubsec}
list explicit expressions for each of the summands in these formulae.
Putting it all together to get actual numbers is a simple computing
exercise.

As for the two leading Lipschitz--Killing curvatures, simple algebra
shows that the results of
this section coincide with those of Corollary~\ref{corscale}.

For notational convenience in what follows, we make the simple scale
transformation
$s=-\ln\sigma$, so that our random field is now denoted by
$f(s,t)$ and has covariance function
\[
C\bigl((s_1,t_1),(s_2,t_2)
\bigr)= e^{N(s_1+s_2)/2} \int_{\mathbb{R}^N} h \bigl((t_1-u)e^{s_1}
\bigr) h \bigl((t_2-u)e^{s_2} \bigr) \,du.
\]
The parameter space of $f$, which we shall denote by $M$, is now of the form
$T\times[s_1,s_2]\subset\real^{N+1}$, for some $-\infty
<s_1<s_2<\infty$.\vadjust{\goodbreak}

\subsection{Stratifying the parameter space}

We start
by stratifying the parameter space $M$ into manifolds of common
dimension. In
particular, we write the $(N+1)$-dimensional parameter space
$
M = T \times[s_1, s_2]
$
as the disjoint union of four types of pieces, arranged in three
strata, according
to their dimension (recall that $T$ has a smooth boundary $\partial
T$):
\begin{eqnarray*}
\partial_{N+1}M &=& {M}^{\circ} = (s_1,
s_2) \times{T}^{\circ
},
\\
\partial_{N}M &=& (s_1, s_2) \times\partial T
\cup \{s_2\} \times{T}^{\circ}  \cup
\{s_1\} \times{T}^{\circ}
\\
&\definedas& \mbox{``side'' $\cup$
``top''  $\cup$ ``bottom''},
\\
\partial_{N-1}M &=& \{s_1\} \times\partial T \cup
\{s_2\} \times \partial T.
\end{eqnarray*}
%

\subsection{What we need to compute}

Ultimately, we need to compute the Lipschitz--Killing curvatures $\lips_i(M)$
of $M$ under the Riemannian metric induced
by the random field $f$. Rather than doing this directly, we shall compute
certain other curvatures, denoted by $\lips^\a_i(M)$,
specifically designed to be simpler to handle on sets of
constant curvature. (We shall soon see that
$M$ has constant curvature, given by~$-\kappa^{-1}$.)

There are simple relationships between the $\lips_i(M)$ and $\lips^\a_i(M)$,
including the following one, which is (10.5.12) of~\cite{TheBook}:
%
\begin{equation}
\label{eqdefcurvaturemeasureonetubes} \lips_i(M) = \sum
_{n=0}^{\lfloor({N+1-i})/{2}\rfloor} \frac{\a^n(i+2n)!}{(4\pi)^{n}n!i!}
\lips^\a_{i+2n}(M).
\end{equation}

To compute the $\lips^\a_i(M)$ themselves, we
first write
%
\begin{equation}
\label{ECsumofLsequn} \lips^\a_i(M) = \sum
_{j=N-1}^{N+1} \lips^\a_i(M;
\partial_j M)
\end{equation}
and then restrict ourselves to the case $\a=-\kappa^{-1}$.
Applying (10.7.10)\setcounter{footnote}{4}\footnote{There is a small, but, for us, rather
significant typo in (10.7.10) of~\cite{TheBook},
in that there is a minus sign missing before the $\kappa$ in $\kappa I^2/2$.}
of~\cite{TheBook} along with the (yet to be proven) fact that $M$ has
constant curvature
$-\kappa^{-1}$ gives that, for $i\leq j$,
%
\begin{eqnarray}
\label{eqcurvaturemeasuresgaussiangeneral}\qquad
&&\lips^{-\kappa^{-1}}_i(M;
\partial_j M)
\nonumber\\[-8pt]\\[-8pt]
&&\qquad = \frac
{1}{(2\pi)^{(j-i)/2}(j-i)!} \int_{\partial_jM } \E \bigl\{
\Tr^{T_t\partial_jM}\bigl({S}_{Z_{j,t}}^{j-i}\bigr)
\indic_{N_tM}(Z_{j,t}) \bigr\} \hauss_j(dt).
\nonumber
\end{eqnarray}
Here, for each $t \in\partial_jM$, $Z_{j,t}$ is a normally distributed
random vector of dimension $N+1-j$ in the space $T_t\partial_jM^{\perp
}$, the
orthogonal complement of $T_t\partial_jM$ in $\real^{N+1}$, $S$~is the
Riemannian shape
operator of $M$, $N_tM$ is the normal cone to $M$ at $t$,
and $\hauss_j$ is the volume measure on $\partial_jM$ corresponding
to the Riemannian metric induced by the random field $f$.

The next step is to compute each element in
(\ref{eqcurvaturemeasuresgaussiangeneral}).

\subsection{The induced Riemannian metric}

The next step is to identify the Riemannian metric
$g$ that the random field induces on $M$. In order to describe this, however,
we need first to choose a family of vector fields generating the
tangent bundle of $M$.

We do this sequentially, starting with the set $T$, ignoring for the
moment the scale
component of the parameter space. For $T$, let $\eta$ be the
(Euclidean) outward unit normal vector
field on $\partial T$, and extend this to a full (Euclidean) $C^2$, orthonormal
tangent vector bundle
$X_1, \ldots, X_N$, on $T$, with $X_N=\eta$. Enlarge this to a vector
bundle on all of $M$ by
adding the vector field $\nu=\partial/\partial s$, the field of
tangent vectors in the scale
direction.

The Riemannian metric $g$ at points $(s,t)\in M$ can now be calculated on
pairs of vectors from the above
vector field using the variances (\ref{ECscale-variancesequn})
and (\ref{ECscale-variances2equn}) and the independence discussed there.
These yield
%
\begin{eqnarray}
\label{ECRM1equn}
g_{(s,t)}(X_{i,(s,t)}, X_{j,(s,t)})
&=&
\gamma e^{2s} \delta_{ij},
\\
\label{ECRM2equn}
g_{(s,t)}(X_{i,(s,t)}, \nu_{(s,t)}) &=& 0,
\\
\label{ECRM3equn}
g_{(s,t)}(\nu_{(s,t)}, \nu_{(s,t)}) &=& \kappa,
\end{eqnarray}
where $\delta_{ij}$ is the Kronecker delta. We note, once again,
that (\ref{ECRM2equn})
implies that the metric $g$ is actually of product form is essential to
all that follows.

Note that the above three equations also
imply that the structure of the normal cones
to $M$ is well described by the initial, Euclidean, choice
of vector fields. There is, of course, no normal cone at points $t\in
\partial_{N+1}M$, the
interior of $M$. For the sets in $\partial_NM$, one of the vector
fields $\eta$ and
$\nu$ describes the
normal geometry. In particular, along the side the normal is an
$s$-dependent multiple of
$\eta$ and along the top and bottom the normals are constant multiples
of $\nu$.
Along $\partial_{N-1}$ all normals are linear combinations of elements
in the vector
fields $\eta$ and $\nu$.

\subsection{The Levi-Civita connection}

The next step involves computing the second fundamental forms in
(\ref{eqcurvaturemeasuresgaussiangeneral}) and so requires
identifying the
Levi-Civita connections. If we think of $M$ as embedded in a smooth
open subset
of $\real^{N+1}$, with the same induced metric,
then we shall write
$\nabla$ for the Levi-Civita connection on $M$ and $\widetilde\nabla
$ for the Levi-Civita connection
on the ambient set.

For our purposes, we need only to know how $\widetilde\nabla$ operates
on vectors normal to~$M$. Furthermore, since the Lie bracket
%
\begin{equation}
\label{ECXinuzeroequn} [X_i,\nu] = 0\qquad \mbox{for all }  1 \leq i
\leq N,
\end{equation}
we need only compute $\widetilde{\nabla}_{X_i}X_j$, $\widetilde
{\nabla}_{\nu}X_i$ and $\widetilde{\nabla}_{\nu}\nu$.

To start, note that a straightforward application of Koszul's
formula [which in the current situation simplifies to
$ 2 g( \nabla_X Y, Z ) = X g( Y,Z ) + \break Y g( X, Z ) - Z g( X, Y )$]
and (\ref{ECRM1equn}) and (\ref{ECRM2equn}) imply, for
$i,j,k=1,\ldots,N$, that
\[
g(\widetilde{\nabla}_{X_i}X_j, X_k) = \gamma
e^{2s} \langle\nabla_{X_i}X_j, X_k
\rangle= g( \nabla_{X_i}X_j, X_k),
\]
where $\nabla$ is the standard Euclidean connection, and
%
\begin{equation}
\label{ECconnection0equn} g(\widetilde{\nabla}_{X_i}X_j,
\nu) = -\frac{\gamma}{2} \,\frac
{\partial}{\partial s}e^{2s}
\delta_{ij} = -\gamma e^{2s} \delta_{ij}.
\end{equation}
Now note that for any tangent vectors $X,Y$ to $T$, $\widetilde\nabla_XY$ is a vector
in $\real^{N+1}$, and so can be written as $\sum_{j=1}^Na_jX_j + b\nu
$ for appropriate
coefficients. This fact, together with the last two equalities and
(\ref{ECRM3equn}), implies
that for any two vector fields $X, Y$ on $T$,
%
\begin{equation}
\label{ECconnection1equn} \widetilde{\nabla}_XY =
\nabla_XY - \kappa^{-1} \gamma e^{2s} \langle
X, Y \rangle\nu,
\end{equation}
giving us the first computation.

As far as $\widetilde\nabla_\nu\nu$ is concerned, note first that
another easy consequence of Koszul's formula is that
%
\begin{equation}
\label{ECconnection2equn} g(\widetilde{\nabla}_{\nu}\nu, \nu) =
\tfrac{1}{2} \nu\bigl( g(\nu,\nu)\bigr) = 0,
\end{equation}
giving us $\widetilde\nabla_\nu\nu\equiv0$.

All that remains is to compute $\widetilde{\nabla}_{\nu}X_i$.
An application of the Weingarten equation
[i.e., the scalar second fundamental form is given by
$
S_\nu(X,Y) =
g(\widetilde\nabla_XY, \nu)
= - g(Y, \widetilde\nabla_X\nu)
$]
and (\ref{ECconnection1equn}) yield
\[
g( \widetilde{\nabla}_{\nu}X_i, X_j) =
\frac{\gamma}{2} \,\frac
{\partial}{\partial s} e^{2s} \delta_{ij} =
\gamma e^{2s} \delta_{ij} = g(X_i,X_j)
\]
and $
g( \widetilde{\nabla}_{\nu}X_i, \nu) = 0.
$
Applying now the torsion freeness of connections (i.e.,
$\nabla_XY - \nabla_YX - [X,Y] = 0$), along with (\ref{ECXinuzeroequn})
to the above, gives
%
\begin{equation}
\label{eqmixedconnection} \widetilde\nabla_{\nu}X_i =
\widetilde\nabla_{X_i}\nu= X_i,
\end{equation}
and we have the last of the three cases we were seeking.

\subsection{The second fundamental forms and the curvature matrix}

With the relevant connections determined, we can now turn to
computing second fundamental forms and curvature matrices along the three
stratifications of $M$.

1. $\partial_{N+1}M$, \textit{the interior
$(s_1, s_2) \times{T}^{\circ}$}:
Since the normal space in the interior is empty, there are no normals, no
second fundamental form and no curvature matrix here.

2. $\partial_NM$: \textit{The side $(s_1,s_2)
\times\partial T$}:
A convenient choice for an orthonormal (in the metric $g$)
basis for the tangent space at any point on the side is
given~by\looseness=-1
\[
\biggl\{\frac{X_1}{\gamma^{1/2}e^{s}}, \ldots, \frac{X_{N-1}}{\gamma^{1/2}e^{s}}, \frac{\nu}{\kappa^{1/2}}
\biggr\}
\]\looseness=0
with outward unit normal vector $ \eta/\gamma^{1/2}e^{s}$.\vadjust{\goodbreak}

The scalar second fundamental forms of interest are therefore
\begin{eqnarray*}
&\displaystyle S_{{\eta}/({\gamma^{1/2}e^{s}})} \biggl(\frac{X_i}{\gamma^{1/2}e^{s}}, \frac{X_j}{\gamma^{1/2}e^{s}}
\biggr),\qquad
S_{{\eta}/({\gamma^{1/2}e^{s}})} \biggl(\frac{X_i}{\gamma^{1/2}e^{s}}, \frac{\nu}{\kappa^{1/2}}
\biggr),&\\
&\displaystyle S_{{\eta}/({\gamma^{1/2}e^{s}})} \biggl( \frac{\nu}{\kappa^{1/2}}, \frac{\nu}{\kappa^{1/2}}
\biggr),&
\end{eqnarray*}
all of which can be computed directly using
Weingarten's equation
and (\ref{ECconnection1equn})--(\ref{eqmixedconnection}).
We summarize them in a curvature matrix, writing
\[
\pmatrix{ C = \gamma^{-1/2}e^{-s} C_t & 0
\cr
0 & 0 },
\]
where $C_t$ is the $(N-1)\times(N-1)$
Euclidean curvature matrix of $\partial T$ in the basis
$\{X_1, \ldots, X_{N-1}, \eta\}$.
Recall for later use that, from Section 7.2 of~\cite{TheBook}, in the
current scenario,
the trace in (\ref{eqcurvaturemeasuresgaussiangeneral}) can be replaced
by
%
\begin{equation}
\label{ECtrace-detequn} (j-i)! \detr_{j-i}C,
\end{equation}
where $\detr_j$ is our usual
sum of determinants of $j \times j$ principal minors.

3. $\partial_NM$, \textit{the bottom $\{s_1\}
\times{T}^{\circ}$ and the top
$\{s_2\} \times{T}^{\circ}$}:
Starting with the bottom $\{s_1\} \times{T}^{\circ}$, the outward
unit normal vector is $-\kappa^{-1/2}\nu$, and a convenient
orthonormal basis is given by
$
e^{-s}\{{X_1}, \ldots, {X_{N}}\}
$.
The curvature matrix is therefore $N \times N$ with entries
\[
S_{-\kappa^{-1/2}\nu} \biggl( \frac{X_i}{\gamma^{1/2}e^s}, \frac{X_j}{\gamma^{1/2}e^s} \biggr)= g
\biggl(\widetilde{\nabla}_{{X_i}/({\gamma^{1/2}e^s})}\bigl(-\kappa^{-1/2}\nu\bigr),
\frac{X_j}{\gamma^{1/2}e^s} \biggr) = \kappa^{-1/2} \delta_{ij}
\]
by Weingarten's equation and (\ref{ECconnection0equn}).

Thus, for the bottom, the outward curvature matrix is $\kappa^{-1/2}
I_{N\times N}$,
while along the top
the same arguments give it as $-\kappa^{-1/2} I_{N\times N}$.

4. $\partial_{N-1}M$: \textit{The edges
$\{s_1\} \times\partial T$ and $\{s_2\} \times\partial T$}:
For the edges, we need to consider the scalar second fundamental form
itself and not just the curvature matrix.
As above, an orthonormal basis for the tangent space is
$\{{X_1}/{\gamma^{1/2}e^{s}}, \ldots,\allowbreak {X_{N-1}}/{\gamma^{1/2}e^{s}},
{\nu}/{\kappa^{1/2}}\}$,
but\vspace*{1pt} now an orthonormal basis for the normal space is
$\{\gamma^{-1/2}e^{-s}\eta,\kappa^{-1/2}\nu\}$.

Applying the Weingarten equation, we need to compute
\[
g \bigl(\widetilde{\nabla}_{\gamma^{-1/2}e^{-s} X_i} \bigl(\gamma^{-1/2}e^{-s}
X_j \bigr), a \gamma^{-1/2}e^{-s} \eta+ b
\kappa^{-1/2}\nu \bigr) 
\]
for arbitrary $a, b$. Applying
(\ref{ECconnection1equn}) gives that this is
$a \gamma^{-1/2} e^{-s} C_{ij,t} + b \kappa^{-1/2} \delta_{ij}
$, 
where, with a minor abuse of notation, we now use $C_t$
to denote the $(N-1)\times(N-1)$ Euclidean curvature matrix of
$\partial M$ at $t$.

\subsection{The curvature tensor}

We now have all the pieces we need to compute the Lipschitz--Killing
curvatures $\lips_j^{-\kappa^{-1}}$,
and could actually proceed to the final computation. However, in
justifying the
formula (\ref{eqcurvaturemeasuresgaussiangeneral}), we used the
fact that
$M$ has constant negative curvature $-\kappa^{-1}$. Now we shall take
a moment to prove this.
Let
$
R(X,Y) = \widetilde{\nabla}_X \widetilde{\nabla}_Y - \widetilde
{\nabla}_Y \widetilde{\nabla}_X - \widetilde{\nabla}_{[X,Y]}$
be the curvature operator.
From previous calculations of the connection
\begin{eqnarray*}
R(X_i,\nu)X_k &=& \widetilde{\nabla}_{X_i}
\widetilde{\nabla}_{\nu
}X_k - \widetilde{
\nabla}_{\nu}\widetilde{\nabla}_{X_i}X_k
\\[-2pt]
&=& \widetilde{\nabla}_{X_i}X_k - \widetilde{
\nabla}_{\nu} \bigl(\nabla_{X_i}X_k -
\delta_{ik} \kappa^{-1} \gamma e^{2s} \nu\bigr)
\\[-2pt]
&=& \nabla_{X_i}X_k - \delta_{ik}
\kappa^{-1} \gamma e^{2s} \nu- \nabla_{X_i}X_k
+ \delta_{ik} \kappa^{-1} \,\frac{\partial}{\partial
s} \gamma
e^{2s} \nu
\\[-2pt]
&=& \delta_{ik} \kappa^{-1}\gamma e^{2s} \nu,
\\[-2pt]
R(X_i, \nu)\nu&=& \widetilde{\nabla}_{X_i} \widetilde{
\nabla}_{\nu
}\nu- \widetilde{\nabla}_{\nu} \widetilde{
\nabla}_{X_i}\nu
\\[-2pt]
&=& - \widetilde{\nabla}_{\nu}X_i
\\[-2pt]
&=& -X_i,
\\[-2pt]
R(X_i,X_j)\nu&=& \widetilde{\nabla}_{X_i}
\widetilde{\nabla }_{X_j}\nu- \widetilde{\nabla}_{X_j}
\widetilde{\nabla}_{X_i}\nu- \widetilde{\nabla}_{[X_i,X_j]}\nu
\\[-2pt]
&=& \widetilde{\nabla}_{X_i}X_j - \widetilde{
\nabla}_{X_j}X_i - [X_i,X_j]
\\[-2pt]
&=& \nabla_{X_i}X_j - \nabla_{X_j}X_i
- [X_i,X_j] - \kappa^{-1}\gamma
e^{2s} (\delta_{ij} - \delta_{ji}) \nu
\\[-2pt]
&=& 0,
\\[-2pt]
R(X_i,X_j)X_k &=& \widetilde{
\nabla}_{X_i}\widetilde{\nabla }_{X_j}X_k -
\widetilde{\nabla}_{X_j}\widetilde{\nabla}_{X_i}X_k
- \widetilde{\nabla}_{[X_i,X_j]}X_k
\\[-2pt]
&=& \widetilde{\nabla}_{X_i} \bigl(\nabla_{X_j}X_k
- \delta_{jk} \kappa^{-1} \gamma e^{2s} \nu\bigr)
- \widetilde{\nabla}_{X_j} \bigl(\nabla_{X_i}X_k
- \delta_{ik} \kappa^{-1}\gamma e^{2s} \nu\bigr)
\\[-2pt]
&&{} - \nabla_{[X_i,X_j]}X_k + \kappa^{-1}\gamma
e^{2s} \bigl\langle[X_i,X_j], X_k
\bigr\rangle\nu
\\[-2pt]
&=& \nabla_{X_i} \nabla_{X_j}X_k -
\delta_{jk} \kappa^{-1}\gamma e^{2s}
X_i - \kappa^{-1} \gamma e^{2s} \langle
\nabla_{X_j}X_k, X_i \rangle\nu
\\[-2pt]
&&{} - \nabla_{X_j} \nabla_{X_i}X_k +
\delta_{ik} \kappa^{-1} \gamma e^{2s}
X_j - \kappa^{-1} \gamma e^{2s} \langle
\nabla_{X_i}X_k, X_j \rangle\nu
\\[-2pt]
&&{} - \nabla_{[X_i,X_j]}X_k + \kappa^{-1}\gamma
e^{2s} \bigl\langle[X_i,X_j], X_k
\bigr\rangle\nu
\\[-2pt]
&=& \nabla_{X_i} \nabla_{X_j}X_k -
\nabla_{X_j} \nabla_{X_i}X_k -
\nabla_{[X_i,X_j]}X_k
\\[-2pt]
&&{} + \kappa^{-1}\gamma e^{2s} \nu \bigl( \bigl
\langle[X_i,X_j], X_k \bigr\rangle- \langle
\nabla_{X_j}X_k, X_i \rangle+ \langle
\nabla_{X_i}X_k, X_j \rangle \bigr)
\\[-2pt]
&&{} - \delta_{jk} \kappa^{-1}\gamma e^{2s}
X_i + \delta_{ik} \kappa^{-1}\gamma
e^{2s} X_j
\\[-2pt]
&=& - \delta_{jk} \kappa^{-1}\gamma e^{2s}
X_i + \delta_{ik} \kappa^{-1}\gamma
e^{2s} X_j,
\end{eqnarray*}
the last equality following from the flatness of $\real^{N+1}$
and torsion freeness.

Since the curvature tensor is given by
$
R(X,Y,Z,W) = g(R(X,Y)Z,W)
$,
it is now easy to use the above calculations to check cases and see that
\[
R(X,Y,Z,W) = -\kappa^{-1}I^2(X,Y,Z,W)/2,\vadjust{\goodbreak}
\]
where $I$ is the identity form given by $I(X,Y)=g(X,Y)$, and with the
usual tensor
product $I^2(X,Y,Z,W)\definedas I(X,Z)I(Y,W)-I(X,W)I(Y,Z)$.

From this follows our claim that $M$ is a space of constant curvature
$-\kappa^{-1}$.\vspace*{-2pt}

\subsection{\texorpdfstring{The Lipschitz--Killing curvatures $\lips^{-\kappa^{-1}}$}
{The Lipschitz--Killing curvatures L -kappa -1}}
\label{scaleLKCSsubsubsec}

With all the preparation done, we can now begin the computation of the
Lipschitz--Killing curvatures $\lips^{-\kappa^{-1}}_j(M;\allowbreak\partial_kM)$, according to their
definition in (\ref{eqcurvaturemeasuresgaussiangeneral}). As in the
discussion of second fundamental forms, we divide the computation into separate
sections, each corresponding to a different stratum in the
stratification of
$M$.

Throughout the following computation we take $Z_1$ and $Z_2$ to be two
independent $N(0,1)$ random variables.
We shall associate the $Z_j$ with the normal vector fields to obtain normal
vectors of the form $Z_1\nu+Z_2\eta$.

We also adopt the notation
%
\begin{eqnarray}
\label{mukequn} \mu_k&\definedas&\E \bigl\{Z_j^k
\indic_{Z_j\geq0} \bigr\} \nonumber\\[-9pt]\\[-9pt]
&=& \cases{\displaystyle  \frac{1}{2}, &\quad if $k=0$,
\vspace*{2pt}\cr
\displaystyle \frac{2^{n-1}n!}{\sqrt{2\pi}}, &\quad if $k=2n+1$ is odd,
\vspace*{2pt}\cr
\displaystyle \frac{(2n-1)(2n-3)\cdots}{2}, &\quad if $k=2n$ is
even.}\nonumber
\end{eqnarray}

1. $\partial_{N+1}M$, \textit{the interior $(s_1, s_2) \times
{T}^{\circ}$}:
Since the second fundamental form is zero in $\partial_{N+1}M$, the
only nonzero
Lipschitz--Killing curvature occurs when $N+1=j=i$ in (\ref
{eqcurvaturemeasuresgaussiangeneral}).
In this case
%
\begin{eqnarray}
\label{ECLKCkappa1equn} \lips_{N+1}^{-\kappa^{-1}} (M;
\partial_{N+1}M ) &=& \hauss_{N+1}\bigl((s_1,s_2)
\times{T}^{\circ}\bigr)
\nonumber
\\[-2pt]
&=& \kappa^{1/2} \gamma^{N/2} \int_{s_1}^{s_2}
\int_{{T}^{\circ}} e^{Ns} \,dt \,ds
\\[-2pt]
&=& \kappa^{1/2} \gamma^{N/2} \frac{e^{Ns_2} - e^{Ns_1}}{N}
\lips_N^E\bigl(T; {T}^{\circ}\bigr),
\nonumber
\end{eqnarray}
where $\lips_N^E(T; {T}^{\circ})$ is computed in the standard
Euclidean sense.


2. $\partial_NM$: \textit{The side $(s_1,s_2) \times\partial T$}:
For this case, all the Lipschitz--Killing curvatures need to be
computed, and we replace the
trace in
(\ref{eqcurvaturemeasuresgaussiangeneral}) by the determinants of
the curvature
matrix as in (\ref{ECtrace-detequn}). Then
%
\begin{eqnarray}
\label{ECLKCkappa2equn} && \lips_j^{-\kappa^{-1}}\bigl(M;
(s_1,s_2) \times\partial T\bigr)
\nonumber\\[-2pt]
&&\qquad = (2\pi)^{-(N-j)/2} \Ee \bigl\{Z_2^{N-j}
\indic_{\{Z_2 > 0\}} \bigr\}
\nonumber\\
&&\qquad\quad{} \times\int_{(s_1,s_2) \times\partial T} \detr_{N-j} \pmatrix{
e^{-s} C_t & 0
\cr
0 & 0 } \,d\hauss_N(s,t)
\nonumber\\[-9pt]\\[-9pt]
&&\qquad = (2\pi)^{-(N-j)/2} \kappa^{1/2}\Ee \bigl
\{Z_2^{N-j}\indic_{\{Z_2 > 0\}
} \bigr\}
\nonumber\\[-2pt]
&&\qquad\quad{} \times\int_{s_1}^{s_2} \int
_{\partial T} \bigl(\gamma^{-1/2}e^{-s}
\bigr)^{N-j} \bigl(\gamma^{1/2}e^s
\bigr)^{N-1} \detr_{N-j}(C_t) \,dt \,ds
\nonumber\\[-2pt]
&&\qquad = \kappa^{1/2}\gamma^{(j-1)/2} \frac
{e^{(j-1)s_2}-e^{(j-1)s_1}}{j-1}
\lips^E_{j-1}(T; \partial T)\nonumber
\end{eqnarray}
in a parallel notation to (\ref{ECLKCkappa1equn}).

3. $\partial_NM$, \textit{the bottom $\{s_1\} \times{T}^{\circ
}$ and the top
$\{s_2\} \times{T}^{\circ}$}:
Beginning with the bottom, $\{s_1\} \times{T}^{\circ}$, for $0 \leq
j \leq N$
we have
%
\begin{eqnarray}
\label{ECLKCkappa3equn} && \lips_j^{-\kappa^{-1}}\bigl(M;
\{s_1\} \times{T}^{\circ}\bigr)
\nonumber\\[-2pt]
&&\qquad = (2\pi)^{-(N-j)/2} \Ee \bigl\{Z_1^{N-j}
\indic_{\{Z_1 > 0\}
} \bigr\}
\nonumber\\[-9pt]\\[-9pt]
&&\qquad\quad{} \times\int_{\{s_1\} \times
{T}^{\circ}} \detr_{N-j} \bigl(
\kappa^{1/2}I_{N \times N} \bigr) \,d\hauss_{N}(t)
\nonumber
\\[-2pt]
&&\qquad = \bigl(2\pi\kappa^{-1}\bigr)^{-(N-j)/2}
\mu_{N-j} \pmatrix{N
\cr
j} \gamma^{N/2} e^{Ns_1}
\lips_N^E\bigl(T; {T}^{\circ}\bigr).\nonumber
\end{eqnarray}
A similar result holds for the top, $\{s_2\} \times{T}^{\circ}$, namely,
%
\begin{eqnarray}
\label{ECLKCkappa4equn}\quad
&& \lips_j^{-\kappa^{-1}}\bigl(M;
\{s_2\} \times{T}^{\circ}\bigr)
\nonumber\\[-9pt]\\[-9pt]
&&\qquad = (-1)^{N-j}\bigl(2\pi\kappa^{-1}
\bigr)^{-(N-j)/2} \mu_{N-j} \pmatrix{N
\cr
j} \gamma^{N/2}
e^{Ns_2}\lips_N^E\bigl(T; {T}^{\circ}
\bigr).\nonumber
\end{eqnarray}

4. $\partial_{N-1}M$: \textit{The edges
$\{s_1\} \times\partial T$ and $\{s_2\} \times\partial T$}:
We start with the top edge, $\{s_2\} \times\partial T$,
\begin{eqnarray*}
&& \lips_j^{-\kappa^{-1}}\bigl(M; \{s_2\} \times
\partial T\bigr)
\\[-2pt]
&&\qquad = (2\pi)^{-(N-1-j)/2} \\[-2pt]
&&\qquad\quad{}\times\int_{\{s_2\} \times\partial T} \Ee \bigl\{
\detr_{N-1-j}\bigl( Z_1 \kappa^{-1/2} I -
Z_2 \gamma^{-1/2}e^{-s} C_t \bigr)
\\[-2pt]
&&\hspace*{161pt}\qquad\quad{} \times\indic_{\{Z_1 > 0\}} \indic_{\{Z_2 > 0\}} \bigr\} \,d
\hauss_{N-1}(s,t).
\end{eqnarray*}
Use now the easily checked expansion that, for $0\leq k\leq n$,
\[
\detr_k (\alpha I_{n\times n}+ A_{n\times n} ) = \sum
_{m=0}^k \alpha^{k-m}
\pmatrix{n-m
\cr
k-m}\detr_{m}(A),
\]
to expand the $\detr$ term in the expectation above and see that
%
\begin{eqnarray}
\label{ECLKCkappa5equn} && \lips_j^{-\kappa^{-1}}\bigl(M;
\{s_2\} \times\partial T\bigr)
\nonumber\\[-2pt]
&&\qquad = (2\pi)^{-(N-1-j)/2} \nonumber\\[-2pt]
&&\qquad\quad{}\times\sum_{m=0}^{N-1-j}
\kappa^{-(N-1-j-m)/2} \pmatrix{N-1-m
\cr
j}
\nonumber\\[-2pt]
&&\hspace*{40pt}\qquad\quad{} \times\Ee \bigl\{Z_1^{N-1-j-m}
\indic_{\{Z_1 > 0\}} \bigr\} \bigl(\gamma^{1/2} e^{s_2}
\bigr)^{N-1} \Ee \bigl\{Z_2^{m}
\indic_{\{Z_2 > 0 \}} \bigr\}
\nonumber\\[-2pt]
&&\hspace*{40pt}\qquad\quad{} \times\int_{\partial T} \detr_{m}\bigl(
\gamma^{-1/2}e^{-s_2}C_t\bigr) \,d
\hauss_{N-1}(t)
\nonumber\\[-9pt]\\[-9pt]
&&\qquad = \bigl(2\pi\kappa^{-1}\bigr)^{-(N-1-j)/2}\nonumber\\[-2pt]
&&\qquad\quad{}\times \sum
_{m=0}^{N-1-j} \bigl(\gamma^{1/2}e^{s_2}
\bigr)^{N-1-m}\pmatrix{N-1-m
\cr
j} \kappa^{m/2}
\nonumber\\[-2pt]
&&\hspace*{40pt}\qquad\quad{} \times\mu_{N-1-j-m} \Ee\bigl\{Z_2^{m}
\indic_{\{
Z_2 > 0 \}}\bigr\} \int_{\partial T} \detr_{m}(C_t)
\,d\hauss_{N-1}(t)
\nonumber\\[-2pt]
&&\qquad = \kappa^{(N-1-j)/2} \bigl(\gamma^{1/2}e^{s_2}
\bigr)^{N-1} \nonumber\sum_{m=0}^{N-1-j}
\pmatrix{N-1-m
\cr
j} \bigl(\gamma^{-1/2}e^{-s_2}
\kappa^{1/2}\bigr)^m\\[-2pt]
&&\hspace*{139pt}\qquad\quad{} \times
\mu_{N-1-j-m} \lips^E_{N-1-m}(T;
\partial T).\nonumber
\end{eqnarray}

A similar argument also works for the bottom edge $\{s_1\} \times
\partial T$,
the only change being that $\indic_{\{Z_1>0\}}$
becomes $\indic_{\{Z_1<0\}}$, giving the final form
%
\begin{eqnarray}
\label{ECLKCkappa6equn} && \lips_j^{-\kappa^{-1}}\bigl(M;
\{s_1\} \times\partial T\bigr)
\nonumber\\[-2pt]
&&\qquad = \kappa^{(N-1-j)/2} \bigl(\gamma^{1/2}e^{s_1}
\bigr)^{N-1} \nonumber\\[-9pt]\\[-9pt]
&&\qquad\quad{}\times\sum_{m=0}^{N-1-j}
\pmatrix{N-1-m
\cr
j} \bigl(\gamma^{-1/2}e^{-s_1}
\kappa^{1/2}\bigr)^m
\nonumber\\[-2pt]
&&\hspace*{41.5pt}\qquad\quad{} \times(-1)^{N-1-j-m} \mu_{N-1-j-m} \lips^E_{N-1-m}(T;
\partial T).
\nonumber
\end{eqnarray}

Collecting (\ref{ECLKCkappa1equn})--(\ref{ECLKCkappa6equn})
now gives us all the $\lips^{-\kappa^{-1}}_i(M;\partial_j M)$, from which,
via (\ref{eqdefcurvaturemeasureonetubes}) and (\ref{ECsumofLsequn}),
we can compute the $\lips_i(M)$, as promised.



\printaddresses

\end{document}